\documentstyle[11pt,amssymb,bm]{article}

\begin{document}

\vskip 45 pt

\centerline{\bf CLASSIFICATION THEOREM ON IRREDUCIBLE
REPRESENTATIONS OF} \centerline{\bf THE $q$-DEFORMED ALGEBRA
$U'_{q}({\rm so}_{n})$}

\bigskip

\centerline{N. Z. Iorgov and A. U. Klimyk}

\medskip

\centerline{Bogolyubov Institute for Theoretical Physics,
Metrologichna str., 03143 Kiev, Ukraine}

\centerline{e-mail: iorgov@bitp.kiev.ua, aklimyk@bitp.kiev.ua}

\bigskip

\begin{abstract}
The aim of this paper is to give a complete classification of
irreducible finite dimensional representations of the nonstandard
$q$-deformation $U'_q({\rm so}_n)$ (which does not coincide with
the Drinfeld--Jimbo
quantum algebra $U_q({\rm so}_n)$) of the universal enveloping
algebra $U({\rm so}_n({\Bbb C}))$ of the Lie algebra ${\rm
so}_n({\Bbb C})$
when $q$ is not a root of
unity. These representations are exhausted by irreducible
representations of the classical type and of the nonclassical
type. Theorem on complete reducibility of finite dimensional
representations of $U'_q({\rm so}_n)$ is proved.
\end{abstract}

\medskip


\bigskip

 \centerline{\sc 1.
Introduction}
\medskip

Quantum orthogonal groups, quantum Lorentz groups and the
corresponding quantum algebras are of special interest for modern
mathematical physics. M. Jimbo [19] and V. Drinfeld [3] defined
$q$-deformations (quantum algebras) $U_q(g)$ for all simple
complex Lie algebras $g$ by means of Cartan subalgebras and root
subspaces (see also [18] and [23]). Reshetikhin, Takhtajan and
Faddeev [32] defined quantum algebras $U_q(g)$ in terms of the
quantum $R$-matrix satisfying the quantum Yang--Baxter equation.
However, these approaches do not give a satisfactory presentation
of the quantum algebra $U_q({\rm so}_n)$ from a viewpoint of some
problems in quantum physics and representation theory. When
considering representations of the quantum algebras $U_q({\rm
so}_{n+1})$ and $U_q({\rm so}_{n,1})$ we are interested in
reducing them onto the quantum subalgebra $U_q({\rm so}_n)$. This
reduction would give an analogue of the Gel'fand--Tsetlin basis
for these representations. However, definitions of quantum
algebras mentioned above do not allow the inclusions $U_q({\rm
so}_{n+1})\supset U_q({\rm so}_n)$ and $U_q({\rm so}_{n,1})\supset
U_q({\rm so}_n)$. To be able to exploit such reductions we have to
consider $q$-deformations of the Lie algebra ${\rm so}_{n+1}({\bf
C})$ defined in terms of the generators
$I_{k,k-1}=E_{k,k-1}-E_{k-1,k}$  (where $E_{is}$ is the matrix
with entries $(E_{is})_{rt}=\delta _{ir} \delta _{st})$ rather
than by means of Cartan subalgebras and root elements. To
construct such deformations we have to deform trilinear relations
for elements $I_{k,k-1}$ instead of Serre's relations (used in the
case of the standard quantized universal enveloping algebras). As
a result, we obtain the associative algebra which will be denoted
as $U'_q({\rm so}_n).$

This $q$-deformation was first constructed in [8]. It permits one
to construct the reductions of $U'_q({\rm so}_{n,1})$ and
$U'_q({\rm so}_{n+1})$ onto $U'_q({\rm so}_n)$. The $q$-deformed
algebra $U'_q({\rm so}_n)$ leads for $n=3$ to the $q$-deformed
algebra $U'_q({\rm so}_3)$ defined by D. Fairlie [4]. The
cyclically symmetric algebra, similar to Fairlie's one, was also
considered somewhat earlier by Odesskii [31].

In the classical case, the imbedding $SO(n)\subset SU(n)$ (and its
infinitesimal analogue) is of great importance for nuclear physics
and in the theory of Riemannian symmetric spaces. It is well known
that in the framework of quantum groups and  Drinfeld--Jimbo
quantum algebras one cannot construct the corresponding embedding.
The algebra $U'_q({\rm so}_n)$ allows to define such an embedding
[29], that is, it is possible to define the embedding $U'_q({\rm
so}_n)\subset U_q({\rm sl}_n)$, where $U_q({\rm sl}_n)$ is the
Drinfeld-Jimbo quantum algebra.

As a disadvantage of the algebra $U'_q({\rm so}_n)$ we have to
mention the difficulties with Hopf algebra structure.
Nevertheless, $U'_q({\rm so}_n)$ turns out to be a coideal in
$U_q({\rm sl}_n)$ (see [29]) and this fact allows us to consider
tensor products of finite dimensional irreducible representations
of $U'_q({\rm so}_n)$ for many interesting cases (see [13]).

The algebra $U'_q({\rm so}_n)$ and their representations are
interesting in many cases. Main directions of interest are the
following:
\medskip

1. The theory of orthogonal polynomials and special functions
(especially, the theory of $q$-orthogonal polynomials and basic
hypergeometric functions). This direction is not good worked out.
Some ideas of such applications can be found in [22].
\smallskip

2. The algebra $U'_q({\rm so}_n)$ (especially its particular case
$U'_q({\rm so}_3)$) is related to the algebra of observables in
2+1 quantum gravity on the Riemmanian surfaces (see papers [2, 5,
28]).
\smallskip

3. A quantum analogue of the Riemannian symmetric space
$SU(n)/SO(n)$ is constructed by means of the algebra $U'_q({\rm
so}_n)$. This construction is fulfilled in the paper [29] (see
also [24]).
\smallskip

4. A $q$-analogue of the theory of harmonic polynomials
($q$-harmonic polynomials on quantum vector space ${\Bbb R}_q^n$)
is constructed by using the algebra $U'_q({\rm so}_n)$. In
particular, a $q$-analogue of different separations of variables
for the $q$-Laplace operator on ${\Bbb R}_q^n$ is given by means
of this algebra and its subalgebras. This theory is contained in
the papers [17] and [30].
\smallskip

5. The algebra $U'_q({\rm so}_n)$ also appears in the theory of
links in the algebraic topology (see [1]).
\smallskip

6. The algebra $U'_q({\rm so}_n)$ is connected with Yangians (see
[26] and references therein).
\smallskip

7. A new quantum analogue of the Brauer algebra is connected with
the algebra $U'_q({\rm so}_n)$ (see [27]).
\medskip

A large class of finite dimensional irreducible representations of
the algebra $U'_q({\rm so}_n)$ were constructed in [8]. The
formulas of action of the generators of $U'_q({\rm so}_n)$ upon
the basis (which is a $q$-analogue of the Gel'fand--Tsetlin basis)
are given there. A proof of these formulas and some their
corrections were given in [6]. However, finite dimensional
irreducible representations described in [6] and [8] are
representations of the classical type. They are $q$-deformations
of the corresponding irreducible representations of the Lie
algebra ${\rm so}_n$, that is, at $q\to 1$ they turn into
representations of ${\rm so}_n$.

The algebra $U'_q({\rm so}_n)$ has other classes of finite
dimensional irreducible representations which have no classical
analogue. These representations are singular at the limit $q\to
1$. They are described in [15]. The description of these
representations for the algebra $U'_q({\rm so}_3)$ is given in
[9]. A classification of irreducible $*$-representations of real
forms of the algebra $U'_q({\rm so}_3)$ is given in [33]. The
representation theory of $U'_q({\rm so}_n)$ when $q$ is a root of
unity is studied in [16].

In this paper we deal with classification of finite dimensional
irreducible representations of the algebra $U'_q({\rm so}_n)$ when
$q$ is not a root of unity. As mentioned above, there were
constructed irreducible representations of the algebra $U'_q({\rm
so}_n)$ belonging to the classical and to the nonclassical types.
However, it was not known that these representations exhaust all
irreducible finite dimensional representations. We started to
study this problem in [21]. We show there that these
representations are determined by the so called highest weights
(which were defined in [21] and differ from highest weights in the
theory of quantized universal enveloping algebras). However, we do
not know a correspondence between known representations of the
classical and nonclassical types and highest weights. In the
present paper we develop an approach to the problem of
classification from other point of view. Namely, we prove that
each irreducible finite dimensional representation of $U'_q({\rm
so}_n)$ belongs to the set of representations of the classical
type or to the set of representations of the nonclassical type,
constructed before. For proving this we use our previous results
on structure of the algebra $U'_q({\rm so}_n)$ (tensor operators,
Wigner--Eckart theorem, etc). We also need the theorem on complete
reducibility of finite dimensional representations of $U'_q({\rm
so}_n)$. This theorem is proved in this paper. Some ideas from the
theory of representations of the Lie algebra ${\rm so}_n({\Bbb
C})$ and its real forms are also used.

Note that the problem of classification of irreducible finite
dimensional representations of $U'_q({\rm so}_n)$ is much more
complicated than in the case of Drinfeld--Jimbo quantum algebras
since in $U'_q({\rm so}_n)$ we do not have an analogue of a Cartan
subalgebra and root elements. The set of all irreducible finite
dimensional representations of $U'_q({\rm so}_n)$ is wider than in
the case of $U_q({\rm so}_n)$.

Everywhere below we assume that $q$ is not a root of unity.

\bigskip

\centerline{\sc 2. The $q$-deformed algebra $U'_q({\rm so}_n)$}
\medskip

The universal enveloping algebra $U({\rm so}_n({\Bbb C}))$ is
generated by the elements $I_{ij}=E_{ij}-E_{ji}$, $i>j$. But in
order to generate the algebra $U({\rm so}_n({\Bbb C}))$, it is
enough to take only the elements $I_{21}$, $I_{32},\cdots
,I_{n,n-1}$. It is a minimal set of elements necessary for
generating $U({\rm so}_n({\Bbb C}))$. These elements satisfy the
relations
$$
I^2_{i,i-1}I_{i+1,i}-2I_{i,i-1}I_{i+1,i}I_{i,i-1} +
I_{i+1,i}I^2_{i,i-1} =-I_{i+1,i},
$$  $$
I_{i,i-1}I^2_{i+1,i}-2I_{i+1,i}I_{i,i-1}I_{i+1,i} +
I^2_{i+1,i}I_{i,i-1} =-I_{i,i-1},
$$   $$
I_{i,i-1}I_{j,j-1}- I_{j,j-1}I_{i,i-1}=0\ \ \ \ {\rm for}\ \ \ \
|i-j|>1.
$$
The following theorem is true for $U({\rm so}_n({\Bbb C}))$ (see
[20]): {\it The enveloping algebra $U({\rm so}_n({\Bbb C}))$ is
isomorphic to the complex associative algebra (with a unit
element) generated by the elements $I_{21}$, $I_{32},\cdots
,I_{n,n-1}$ satisfying the above relations.}

We make a $q$-deformation of these relations by fulfilling the
deformation of the integer 2 as $ 2\to
[2]_q:=(q^2-q^{-2})/(q-q^{-1})=q+q^{-1}$. As a result, we obtain
the relations
$$
I^2_{i,i-1}I_{i+1,i}-(q+q^{-1})I_{i,i-1}I_{i+1,i}I_{i,i-1} +
I_{i+1,i}I^2_{i,i-1} =-I_{i+1,i}, \eqno (1)
$$   $$
I_{i,i-1}I^2_{i+1,i}-(q+q^{-1})I_{i+1,i}I_{i,i-1}I_{i+1,i} +
I^2_{i+1,i}I_{i,i-1} =-I_{i,i-1}, \eqno (2)
$$  $$
I_{i,i-1}I_{j,j-1}- I_{j,j-1}I_{i,i-1}=0\ \ \ \ {\rm for}\ \ \ \
|i-j|>1. \eqno (3)
$$
The $q$-deformed algebra $U'_q({\rm so}_n)$ is defined as the
complex unital (that is, with a unit element) associative algebra
generated by elements $I_{21}$, $I_{32},\cdots ,I_{n,n-1}$
satisfying relations (1)--(3). It is a $q$-deformation of the
universal enveloping algebra $U({\rm so}_n({\Bbb C}))$, different
from the Drinfeld--Jimbo quantized universal enveloping algebra
$U_q({\rm so}_n)$. For this algebra the inclusions $U'_q({\rm
so}_n) \supset U'_q({\rm so}_{n-1})$ and $U_q({\rm sl}_{n})\supset
U'_q({\rm so}_{n})$ are constructed, where $U_q({\rm sl}_{n})$ is
the well known Drinfeld--Jimbo quantum algebra (see Introduction).

An analogue of the skew-symmetric matrices $I_{ij}=E_{ij}-E_{ji}$,
$i>j$, constituting a basis of the Lie algebra ${\rm so}_n({\Bbb
C})$, can be introduced into $U'_q({\rm so}_n)$ (see [7] and
[30]). For $k>l+1$ they are defined recursively by the formulas
$$
 I_{kl}:=
[I_{l+1,l},I_{k,l+1}]_{q}\equiv q^{ 1/2}I_{l+1,l}I_{k,l+1}- q^{-
1/2}I_{k,l+1}I_{l+1,l},
$$
The elements $I_{kl}$, $k>l$, satisfy the commutation relations
$$
[I_{lr},I_{kl}]_q=I_{kr},\ \ [I_{kl},I_{kr}]_q=I_{lr},\ \
[I_{kr},I_{lr}]_q=I_{kl} \ \ \ {\rm for}\ \ \  k>l>r, \eqno (4)
$$   $$
[I_{kl},I_{sr}]=0\ \ \ \ {\rm for}\ \ \ k>l>s>r\ \ {\rm and}\ \
k>s>r>l, \eqno (5)
$$    $$
[I_{kl},I_{sr}]_q=(q-q^{-1}) (I_{lr}I_{ks}-I_{kr}I_{sl}) \ \ \
{\rm for}\ \ \ k>s>l>r. \eqno (6)
$$
For $q=1$ they coincide with the corresponding commutation
relations for the Lie algebra ${\rm so}_n({\Bbb C})$.

The algebra $U'_q({\rm so}_n)$ can be also defined as a unital
associative algebra generated by $I_{kl}$, $1\le l<k\le n$,
satisfying the relations (4)--(6). In fact, the relations (4)--(6)
can be reduced to the relations (1)--(3) for $I_{21}$,
$I_{32},\cdots ,I_{n,n-1}$.

The Poincar\'e--Birkhoff--Witt theorem for the algebra $U'_q({\rm
so}_n)$ can be formulated as follows (a proof of this theorem is
given in [16]): {\it The elements
$$
{I_{21}}^{m_{21}}{I_{31}}^{m_{31}}\cdots {I_{n1}}^{m_{n1}}
{I_{32}}^{m_{32}} {I_{42}}^{m_{42}} \cdots {I_{n2}}^{m_{n2}}
\cdots {I_{n,n-1}}^{m_{n,n-1}},\ \ \ \  m_{ij}=0,1,2, \cdots ,
$$
form a basis of the algebra $U'_q({\rm so}_n)$.}

In $U'_q({\rm so}_n)$ the commutative subalgebra ${\cal A}$
generated by the elements $I_{21}, I_{43}, I_{65},\cdots
,I_{n-1,n-2}$ (or $I_{n,n-1}$) can be separated. So, this
subalgebra is generated by $\lfloor n/2 \rfloor$ elements, where
$\lfloor n/2 \rfloor$ is an integral part of the number $n/2$.
However, there exist no root elements in the algebra $U'_q({\rm
so}_n)$ with respect to this commutative subalgebra. This leads to
the fact that properties of $U'_q({\rm so}_n)$ are not similar to
those of the Drinfeld--Jimbo algebra $U_q({\rm so}_n)$.

\bigskip

\centerline{\sc 3. Irreducible representations of the classical
and nonclassical types}
 \medskip

 \noindent
In this section we give known facts on irreducible representations
of $U'_q({\rm so}_n)$, which will be used below. The corresponding
references are given in Introduction.

Two types of irreducible finite dimensional representations are
known for $U'_q({\rm so}_n)$:
\medskip

(a) representations of the classical type;

(b) representations of the nonclassical type.
\medskip

Known irreducible representations of the classical type are
$q$-deformations of the irreducible finite dimensional
representations of the Lie algebra ${\rm so}_n$. There is a
one-to-one correspondence between these irreducible
representations of the algebra $U'_q({\rm so}_n)$ and irreducible
finite dimensional representations of the Lie algebra ${\rm
so}_n$. Moreover, formulas for representations of the classical
type of $U'_q({\rm so}_n)$ turn into the corresponding formulas
for the representations of Lie algebra ${\rm so}_n$ at $q\to 1$.

There exists no classical analogue for representations of the
nonclassical type: representation operators $T(a)$, $a\in
U'_q({\rm so}_n)$, have singularities at $q=1$.

Let us describe known irreducible finite dimensional
representations of the algebras $U'_q({\rm so}_{n})$, $n \ge 3$,
which belong to the classical type. As in the classical case, they
are given by sets ${\bf m}_{n}$ of $\lfloor {n/2}\rfloor$ numbers
$m_{1,n}, m_{2,n},..., m_{\left \lfloor {n/2}\right \rfloor ,n}$
(here $\lfloor {n/2}\rfloor$ denotes the integral part of ${n/2}$)
which are all integral or all half-integral and satisfy the
dominance conditions
$$
m_{1,2k+1}\ge m_{2,2k+1}\ge ... \ge m_{k,2k+1}\ge 0 ,\ \ \
m_{1,2k}\ge m_{2,2k}\ge ... \ge m_{k-1,2k}\ge |m_{k,2k}|
$$
for $n=2k+1$ and $n=2k$, respectively. These representations are
denoted by $T_{{\bf m}_n}$. We take a $q$-analogue of the
Gel'fand--Tsetlin basis in the representation space, which is
obtained by successive reduction of the representation $T_{{\bf
m}_n}$ to the subalgebras $U'_q({\rm so}_{n-1})$, $U'_q({\rm
so}_{n-2})$, $\cdots$, $U'_q({\rm so}_3)$, $U'_q({\rm
so}_2):=U({\rm so}_2)$. As in the classical case, its elements are
labelled by the Gel'fand--Tsetlin tableaux
$$
\{ \alpha_n\} \equiv \left\{ \matrix{ {\bf m}_{n} \cr {\bf
m}_{n-1} \cr \dots \cr {\bf m}_{2}  }
 \right\}
\equiv \{ {\bf m}_{n},\alpha_{n-1}\}\equiv \{{\bf m}_{n} , {\bf
m}_{n-1} ,\alpha_{n-2}\} , \eqno(7)
$$
where, as in the non-deformed case, the components of ${\bf
m}_{s}$ and ${\bf m}_{s-1}$ satisfy the "betweenness" conditions
$$
m_{1,2p+1}\ge m_{1,2p}\ge m_{2,2p+1} \ge m_{2,2p} \ge ...
\ge m_{p,2p+1} \ge m_{p,2p} \ge -m_{p,2p+1}  ,
$$
$$
m_{1,2p}\ge m_{1,2p-1}\ge m_{2,2p} \ge m_{2,2p-1} \ge ... \ge
m_{p-1,2p-1} \ge \vert m_{p,2p} \vert .
$$
Sometimes, the basis elements, defined by a tableau $\{\alpha_{n}
\}$, are denoted as $\vert \alpha_{n-1} \rangle $ or as $\vert
{\bf m}_{n-1},\alpha_{n-2} \rangle $, that is, we shall omit the
first row ${\bf m}_{n}$ in a tableau.

It is convenient to introduce the so-called $l$-coordinates
$$
l_{j,2p+1}=m_{j,2p+1}+p-j+1,  \qquad
l_{j,2p}=m_{j,2p}+p-j ,
$$
for the numbers $m_{i,k}$. The operator $T_{{\bf
m}_n}(I_{2p+1,2p})$ of the representation $T_{{\bf m}_n}$ of
$U'_q({\rm so}_{n})$ acts upon Gel'fand--Tsetlin basis elements,
labelled by (7), as
$$
T_{{\bf m}_n}(I_{2p+1,2p}) | \alpha_n\rangle = \sum^p_{j=1} \frac{
A^j_{2p}(\alpha_{n})} {a(l_{j,2p}) }
            \vert (\alpha_n)^{+j}_{2p}\rangle -
\sum^p_{j=1} \frac{A^j_{2p}((\alpha_n)^{-j}_{2p})} {a(l_{j,2p}-1)}
|(\alpha_{n})^{-j}_{2p}\rangle \eqno(8)
$$
and the operator $T_{{\bf m}_n}(I_{2p,2p-1})$ acts as
$$
T_{{\bf m}_n}(I_{2p,2p-1})\vert \alpha_n\rangle= \sum^{p-1}_{j=1}
\frac{B^j_{2p-1}(\alpha_{n})} {b( l_{j,2p-1})[l_{j,2p-1}]} \vert
(\alpha_{n})^{+j}_{2p-1} \rangle
$$
$$
-\sum^{p-1}_{j=1}\frac {B^j_{2p-1}((\alpha_n)^{-j}_{2p-1})} {b(
l_{j,2p-1}-1)[l_{j,2p-1}-1]} \vert (\alpha_{n})^{-j}_{2p-1}\rangle
+ {\rm i}\, C_{2p-1}(\alpha_n) \vert \alpha_{n} \rangle . \eqno(9)
$$
In these formulas, $(\alpha_{n})^{\pm j}_{s}$ means the tableau
(7) in which $j$-th component $m_{j,s}$ in ${\bf m}_s$ is replaced
by $m_{j,s}\pm 1$, respectively. The coefficients $A^j_{2p},  $
$B^j_{2p-1},$ $C_{2p-1}$, $a$ and $b$ in (8) and (9) are given by
the expressions
$$
A^j_{2p}(\alpha_n) = \left( \frac{\prod_{i=1}^p
[l_{i,2p+1}+l_{j,2p}] [l_{i,2p+1}-l_{j,2p}-1] \prod_{i=1}^{p-1}
[l_{i,2p-1}+l_{j,2p}] [l_{i,2p-1}-l_{j,2p}-1]} {\prod_{i\ne j}^p
[l_{i,2p}+l_{j,2p}][l_{i,2p}-l_{j,2p}]
[l_{i,2p}+l_{j,2p}+1][l_{i,2p}-l_{j,2p}-1]} \right)^{1/2},  \eqno
(10)
$$
$$
B^j_{2p-1}(\alpha_n)=\left( \frac{\prod_{i=1}^p
[l_{i,2p}+l_{j,2p-1}] [l_{i,2p}-l_{j,2p-1}] \prod_{i=1}^{p-1}
[l_{i,2p-2}+l_{j,2p-1}] [l_{i,2p-2}-l_{j,2p-1}]} {\prod_{i\ne
j}^{p-1} [l_{i,2p-1}{+}l_{j,2p-1}][l_{i,2p-1}{-}l_{j,2p{-}1}]
[l_{i,2p-1}{+}l_{j,2p-1}{-}1][l_{i,2p-1}{-}l_{j,2p-1}{-}1]}
\right) ^{1/2} , \eqno(11)
$$
$$
C_{2p-1}(\alpha_n) =\frac{ \prod_{s=1}^p [ l_{s,2p} ]
\prod_{s=1}^{p-1} [ l_{s,2p-2} ]} {\prod_{s=1}^{p-1} [l_{s,2p-1}]
[l_{s,2p-1} - 1] } ,   \eqno(12)
$$
$$
a(l_{j,2p})=\{ (q^{l_{j,2p}+1}+q^{-l_{j,2p}-1})
(q^{l_{j,2p}}+q^{-l_{j,2p}})\} ^{1/2},\ \ \ b(l_{j,2p-1})=(
[2l_{j,2p-1}+1][2l_{j,2p-1}-1])^{1/2} .
$$
Numbers in square brackets in formulas (9)--(12) mean $q$-numbers
defined by
$$
[a]\equiv [a]_q := \frac {q^a-q^{-a}}{q-q^{-1}}.
$$
It is seen from formula (12) that the coefficient $C_{2p-1}$
vanishes if $m_{p,2p}\equiv l_{p,2p}=0$.

The following assertion is well-known [8]: {\it The
representations $T_{{\bf m}_n}$ are irreducible. The
representations $T_{{\bf m}_n}$ and $T_{{\bf m}'_n}$ are pairwise
nonequivalent for ${\bf m}_n \ne {\bf m}'_n$.}

Irreducible finite dimensional representations of the nonclassical
type are given by sets $\epsilon := (\epsilon _2,\epsilon
_3,\cdots , \epsilon _n)$, $\epsilon _i=\pm 1$, and by sets ${\bf
m}_{n}$ consisting of $\lfloor {n/2}\rfloor$ {\it half-integral}
(but not integral) numbers $m_{1,n}, m_{2,n}, \cdots$, $m_{\lfloor
n/2\rfloor ,n}$ that satisfy the dominance conditions
$$
m_{1,n}\ge m_{2,n}\ge ... \ge m_{\lfloor n/2\rfloor,n}\ge 1/2.
\eqno (13)
$$
These representations are denoted by $T_{\epsilon,{\bf m}_n}$.

For a basis in the representation space, we use an analogue of the
basis of the previous case. Its elements are labelled by tableaux
(7), where the components of ${\bf m}_{s}$ and ${\bf m}_{s-1}$
satisfy the "betweenness" conditions
$$
m_{1,2p+1}\ge m_{1,2p}\ge m_{2,2p+1} \ge m_{2,2p} \ge ...
\ge m_{p,2p+1} \ge m_{p,2p} \ge 1/2  ,
$$
$$
m_{1,2p}\ge m_{1,2p-1}\ge m_{2,2p} \ge m_{2,2p-1} \ge ...
\ge m_{p-1,2p-1} \ge m_{p,2p}  .
$$
The corresponding basis elements are denoted by the same symbols
as in the previous case. The $l$-coordinates for $m_{j,s}$ are
introduced by the same formulas as before.

The operator $T_{\epsilon,{\bf m}_n}(I_{2p+1,2p})$ of the
representation $T_{\epsilon,{\bf m}_n}$ of $U'_q({\rm so}_{n})$
acts upon the basis elements $|\alpha_n\rangle$ by the formulas
$$
T_{\epsilon,{\bf m}_n}(I_{2p+1,2p}) | \alpha_{n}\rangle =
\delta_{m_{p,2p},1/2}\frac{\epsilon_{2p+1}}{q^{1/2}-q^{-1/2}}
D_{2p}(\alpha_n)|\alpha_{n}\rangle+
$$
$$
+\sum^p_{j=1} \frac{ A^j_{2p}(\alpha_n)} {a'(l_{j,2p})}
            \vert (\alpha_{n})^{+j}_{2p}\rangle -
\sum^p_{j=1} \frac{A^j_{2p}((\alpha_n)^{-j}_{2p})}
{a'(l_{j,2p}-1)} |(\alpha_{n})^{-j}_{2p}\rangle , \eqno(14)
$$
where the summation in the last sum must be from 1 to $p-1$ if
$m_{p,2p}=1/2$, and the operator $T_{{\bf m}_n}(I_{2p,2p-1})$
acts as
$$
T_{\epsilon ,{\bf m}_n}(I_{2p,2p-1})\vert \alpha_{n}\rangle=
\sum^{p-1}_{j=1} \frac{B^j_{2p-1}(\alpha_n)} {b(l_{j,2p-1})
[l_{j,2p-1}]_+} \vert (\alpha_{n})^{+j}_{2p-1} \rangle -
$$
$$
-\sum^{p-1}_{j=1}\frac {B^j_{2p-1}((\alpha_n)^{-j}_{2p-1})}
{b(l_{j,2p-1}-1)[l_{j,2p-1}-1]_+} \vert
(\alpha_{n})^{-j}_{2p-1}\rangle + \epsilon _{2p} {\hat
C}_{2p-1}(\alpha_n) \vert \alpha_{n} \rangle , \eqno(15)
$$
where
$$
[a]_+=(q^a+q^{-a})/(q-q^{-1}).
$$
As before, $(\alpha_n)^{\pm j}_{s}$ means the tableau (7) in which
$j$-th component $m_{j,s}$ in ${\bf m}_s$ is replaced by
$m_{j,s}\pm 1$, respectively. The expressions for $A^j_{2p}$,
$B^j_{2p-1}$ and $b$ are given by the same formulas as in (8) and
(9),
$$
a'(l_{j,2p})=\{ (q^{l_{j,2p}+1}-q^{-l_{j,2p}-1})
(q^{l_{j,2p}}-q^{-l_{j,2p}})\} ^{1/2},
$$
$$
{\hat C}_{2p-1}(\alpha_n) = { \prod_{s=1}^p [ l_{s,2p} ]_+
\prod_{s=1}^{p-1} [ l_{s,2p-2} ]_+ \over
   \prod_{s=1}^{p-1} [l_{s,2p-1}]_+ [l_{s,2p-1} - 1]_+ } ,\ \ \ \ \
D_{2p} (\alpha_n)= \frac{\prod_{i=1}^p [l_{i,2p+1}-\frac 12 ]
\prod_{i=1}^{p-1} [l_{i,2p-1}-\frac 12 ] } {\prod_{i=1}^{p-1}
[l_{i,2p}+\frac 12 ] [l_{i,2p}-\frac 12 ] } .  \eqno (16)
$$

The following assertion is true (see [15]): {\it The
representations $T_{\epsilon ,{\bf m}_n}$ are irreducible. The
representations $T_{\epsilon ,{\bf m}_n}$ and $T_{\epsilon ',{\bf
m}'_n}$ are pairwise nonequivalent for $(\epsilon ,{\bf m}_n)\ne
(\epsilon ',{\bf m}'_n)$. For any admissible $(\epsilon ,{\bf
m}_n)$ and ${\bf m}'_n$ the representations $T_{\epsilon ,{\bf
m}_n}$ and $T_{{\bf m}'_n}$ are pairwise nonequivalent.}
\medskip

{\it Remark.} As in the case of irreducible representations of the
Lie algebra ${\rm so}_n$, it follows from the explicit description
of irreducible representations $T_{{\bf m}_n}$ and $T_{\epsilon
,{\bf m}_n}$ of $U'_q({\rm so}_{n})$ that the restriction of
$T_{{\bf m}_n}$ onto the subalgebra $U'_q({\rm so}_{n-1})$
decomposes into a direct sum of irreducible representations of
this subalgebra belonging to the classical type and the
restriction of $T_{\epsilon ,{\bf m}_n}$ onto $U'_q({\rm
so}_{n-1})$ decomposes into a direct sum of irreducible
representations belonging to the nonclassical type. Formulas for
the representations determine explicitly these decompositions.

\bigskip

\centerline{\sc 4. Vector operators and Wigner--Eckart theorem}
 \medskip

In this section we define vector operators for irreducible
representations of $U'_q({\rm so}_n)$ and give the Wigner--Eckart
theorem for them. This information will be used under proving our
main results.

The algebra $U'_q({\rm so}_{n})$ is not a Hopf algebra. For this
reason, we cannot define a tensor products of its representations.
However, $U'_q({\rm so}_{n})$ can be embedded into the Hopf
algebra $U_q({\rm sl}_{n})$ (see [29] and [30]). Using this
embedding, a tensor product of the irreducible representations
$T_1$ and $T$ of $U'_q({\rm so}_{n})$ is determined, where $T_1$
is a vector representation (that is, a representation of the
classical type characterized by the numbers $(1,0,\cdots ,0)$) and
$T$ is an arbitrary irreducible finite dimensional representation
[13]. The decomposition of this tensor product into irreducible
constituents is given by the formulas as in the classical case if
the representation $T$ belongs to the classical type (that is, the
decomposition of $T_1\otimes T_{{\bf m}_n}$ contains the
irreducible representations of the classical type characterized by
${\bf m}^{+j}_n$, ${\bf m}^{-j}_n$, $j=1,2,\cdots ,\lfloor
n/2\rfloor$, and also the representation $T_{{\bf m}_n}$ if
$n=2k+1$ and $m_{k,2k+1}\ne 0$). For the representations
$T=T_{\epsilon ,{\bf m}_n}$ of the nonclassical type we have
$$
T_1\otimes T_{\epsilon ,{\bf m}_n}=\bigoplus _{{\bf m}'_n\in
S_\epsilon ({\bf m}_n)} T_{\epsilon ,{\bf m}'_n},
$$
where
$$
S_\epsilon ({\bf m}_{2p+1})=\bigcup _{j=1}^p \{ T_{\epsilon ,{\bf
m}^{+j}_{2p+1}} \} \cup \bigcup _{j=1}^p \{ T_{\epsilon ,{\bf
m}^{-j}_{2p+1}}\} \cup \{ T_{\epsilon ,{\bf m}_{2p+1}}\} ,\ \ \
S_\epsilon ({\bf m}_{2p})=\bigcup _{j=1}^p \{ T_{\epsilon ,{\bf
m}^{+j}_{2p}} \} \cup \bigcup _{j=1}^p \{ T_{\epsilon ,{\bf
m}^{-j}_{2p}}\} .
$$
As before, ${\bf m}^{\pm j}_n$ is the set of numbers ${\bf m}_n$
with $m_{jn}$ replaced by $m_{jn}\pm 1$, respectively. Note that
each representation $T_{{\bf m}'_{n}}$ and each representation
$T_{\epsilon ,{\bf m}'_{n}}$ for which ${\bf m}'_{n}$ does not
satisfy the dominance conditions must be omitted. Proofs of these
decompositions can be found in [14]. As in the case of quantized
universal enveloping algebras (see [23], Chapter 7),
decompositions of the above tensor products are fulfilled by means
of matrices whose entries are called {\it Clebsch--Gordan
coefficients}.

Let us define a vector operator (it is a set of $n$ operators)
which transforms under the vector representation of the algebra
$U'_q({\rm so}_{n})$. This operator acts on a linear space ${\cal
H}$ on which some representation $T$ of $U'_q({\rm so}_{n})$ acts.
We shall consider only the case when ${\cal H}$ is a finite
dimensional space. We also suppose that ${\cal H}$ decomposes into
a direct sum of irreducible invariant (with respect to $U'_q({\rm
so}_{n})$) subspaces, where only irreducible representations of
the classical type or only irreducible representations of the
nonclassical type are realized. This assumption is explained by
the fact that a vector operator cannot map a subspace on which an
irreducible representation of the classical type is realized into
a subspace on which a representation of the nonclassical type is
realized, or vise versa.

The set $A_r$, $r=1,2,\cdots ,n$, of operators on ${\cal H}$ is
called a {\it vector operator} for the algebra $U'_q({\rm
so}_{n})$ if
$$
[A_{j-1},T(I_{j,j-1})]_q=A_j, \ \ \
[T(I_{j,j-1}),A_j]_q=A_{j-1},
$$    $$
[T(I_{j,j-1}), A_k]_q=0, \ \ \ k\ne j, j-1,
$$
where $[X,Y]_q\equiv q^{1/2}XY-q^{-1/2}YX$ and $T$ is a fixed
representation of $U'_q({\rm so}_n)$ acting on ${\cal H}$.

We represent the space ${\cal H}$ as a direct sum of irreducible
invariant (with respect to $U'_q({\rm so}_n)$) subspaces
$$
{\cal H}=\bigoplus _{\epsilon,{\bf m}_{n},i} {\cal
V}_{\epsilon,{\bf m}_{n},i},
$$
where ${\cal V}_{\epsilon,{\bf m}_{n},i}$ is a subspace, on which
an irreducible representation of $U'_q({\rm so}_{n})$
characterized by $\epsilon$ and ${\bf m}_{n}$ is realized, and $i$
separates multiple irreducible representations of $U'_q({\rm
so}_{n})$ in the decomposition. If irreducible representations
belong to the classical type, then $\epsilon$ must be omitted.

We take a Gel'fand--Tsetlin basis in each subspace ${\cal
V}_{\epsilon,{\bf m}_{n},i}$ and denote these basis vectors by
$|\epsilon,{\bf m}_{n},i,\alpha \rangle$, where $\alpha\equiv
\alpha_{n-1}$ are the corresponding Gel'fand--Tsetlin tableaux.
Then the subspaces
$$
{\cal V}^\alpha_{\epsilon,{\bf m}_{n}}=\bigoplus _i {\Bbb C}
|\epsilon,{\bf m}_{n},i,\alpha \rangle
$$
can be defined.

The Wigner--Eckart theorem for vector operators $\{ A_j\}$ (proved
in [14]) states that the matrix elements of $A_j$ are of the form
$$
\langle\epsilon',{\bf m}'_{n},i',\alpha' |A_j| \epsilon,{\bf
m}_{n},i,\alpha \rangle= C^{\epsilon',{\bf
m}'_{n},\alpha'}_{j;\epsilon,{\bf m}_{n},\alpha}
\langle\epsilon',{\bf m}'_{n},i' \Vert A \Vert \epsilon,{\bf
m}_{n},i \rangle , \eqno (17)
$$
where $C^{\epsilon',{\bf m}'_{n-1},\alpha'}_{j;\epsilon,{\bf
m}_{n-1},\alpha}$ are Clebsch--Gordan coefficients of the tensor
product $T_1\otimes T_{\epsilon,{\bf m}_{n}}$ (these coefficients
are given in an explicit form in [14]), and $\langle\epsilon',{\bf
m}'_{n-1},i' \Vert A \Vert \epsilon,{\bf m}_{n-1},i \rangle $ are
called {\it reduced matrix elements} of the vector operator $\{
A_j\}$. These reduced matrix elements depend only on numbers
characterizing the representations and on the indices separating
multiple representations, and are independent of basis elements of
irreducible invariant subspaces. They are also independent of the
number $j$ of the operator $ A_j$. In the above formulas,
$\epsilon$ must be omitted if we deal only with representations of
the classical type.

Due to the formulas for decompositions of the tensor products
$T_1\otimes T_{{\bf m}_{n}}$ and $T_1\otimes T_{\epsilon,{\bf
m}_{n}}$ we find that matrix elements $\langle\epsilon',{\bf
m}'_{n},i',\alpha' |A_j| \epsilon,{\bf m}_{n},i,\alpha \rangle$
can be non-vanishing only if $\epsilon'= \epsilon$ and also ${\bf
m}'_{n}= {\bf m}^{\pm s}_{n}$ or ${\bf m}'_{n}={\bf m}_{n}$ (since
only for these cases the corresponding Clebsch--Gordan
coefficients can be non-vanishing). Due to the above formulas for
decompositions of tensor products of representations, a vector
operator cannot map a subspace of an irreducible representation of
the classical type (of the nonclassical type) into subspaces on
which irreducible representations of the nonclassical type (of the
classical type) are realized. Therefore, in matrix elements (17)
both indices $\epsilon$ and $\epsilon'$ exist or both are absent.

We can define the operators
$$
A^{{\bf m}_{n}}_{{\bf m}_{n}}: {\cal V}^\alpha_{\epsilon,{\bf
m}_{n}}\to {\cal V}^{\alpha}_{\epsilon,{\bf m}_{n}},\ \ \ \ \
A^{{\bf m}^{+j}_{n}}_{{\bf m}_{n}} : {\cal
V}^\alpha_{\epsilon,{\bf m}_{n}}\to {\cal
V}^{\alpha'}_{\epsilon,{\bf m}^{+j}_{n}},\ \ \ \ \  A^{{\bf
m}^{-j}_{n}}_{{\bf m}_{n}} : {\cal V}^\alpha_{\epsilon,{\bf
m}_{n}}\to {\cal V}^{\alpha'}_{\epsilon,{\bf m}^{-j}_{n}}
$$
which have matrix elements coinciding with reduced matrix elements
of the tensor operator $\{ A_j\}$:
$$
\langle \epsilon,{\bf m}_{n},i',\alpha |A^{{\bf m}_{n}}_{{\bf
m}_{n}}| \epsilon,{\bf m}_{n},i,\alpha \rangle= \langle
{\epsilon,\bf m}_{n},i' \Vert A \Vert \epsilon,{\bf m}_{n},i
\rangle ,
$$   $$
\langle \epsilon,{\bf m}^{+j}_{n},i',\alpha' |A^{{\bf
m}^{+j}_{n}}_{{\bf m}_{n}}| \epsilon,{\bf m}_{n},i,\alpha \rangle=
\langle \epsilon,{\bf m}^{+j}_{n},i' \Vert A \Vert \epsilon,{\bf
m}_{n},i \rangle ,
$$    $$
\langle \epsilon,{\bf m}^{-j}_{n},i',\alpha' |A^{{\bf
m}^{-j}_{n}}_{{\bf m}_{n}}| \epsilon,{\bf m}_{n},i,\alpha \rangle=
\langle \epsilon,{\bf m}^{-j}_{n},i' \Vert A \Vert \epsilon,{\bf
m}_{n},i \rangle .
$$
(The symbol $\epsilon$ must be omitted in these formulas if
necessary.) It follows from the Wigner--Eckart theorem that for
any irreducible representation $T_{\epsilon,{\bf m}_n}$ contained
in the representation $T$, these operators satisfy the following
relations
$$
T_{\epsilon,{\bf m}_n}(a)A^{{\bf m}_{n}}_{{\bf m}_{n}}=A^{{\bf
m}_{n}}_{{\bf m}_{n}} T_{\epsilon,{\bf m}_n}(a), \ \ \ \ a\in
U'_q({\rm so}_n),
$$   $$
T_{\epsilon,{\bf m}_n}(a)A^{{\bf m}_n}_{{\bf m}^{\mp
j}_{n}}A^{{\bf m}^{\pm j}_{n}}_{{\bf m}_{n}}= A^{{\bf
m}_{n}}_{{\bf m}^{\mp j}_{n}}A^{{\bf m}^{\pm j}_{n}}_{{\bf
m}_{n}}T_{\epsilon,{\bf m}_n}(a) ,\ \ \ \ a\in U'_q({\rm so}_n),
$$
where $A^{{\bf m}_n}_{{\bf m}^{\mp j}_{n}}A^{{\bf m}^{\pm
j}_{n}}_{{\bf m}_{n}}$ is considered as operators from ${\cal
V}^\alpha_{\epsilon,{\bf m}_{n}}$ into ${\cal
V}^\alpha_{\epsilon,{\bf m}_{n}}$.
 \medskip

{\bf Proposition 1.} {\it Let $\xi\in {\cal H}$ belongs to a
subspace ${\cal H}_{{\bf m}_n}$ of the irreducible representation
$T_{{\bf m}_n}$ of $U'_q({\rm so}_{n})$. Then $A^{{\bf m}^{+
j}_{n}}_{{\bf m}_{n}}\xi$ and $A^{{\bf m}^{-j}_{n}}_{{\bf
m}_{n}}\xi$ belong to some subspaces ${\cal H}_{{\bf m}^{+j}_n}$
and ${\cal H}_{{\bf m}^{-j}_n}$ of ${\cal H}$, on which the
irreducible representations $T_{{\bf m}^{+j}_n}$ and $T_{{\bf
m}^{-j}_n}$ of $U'_q({\rm so}_{n})$ are realized, respectively.
All the vectors $A^{{\bf m}^{\pm j}_{n}}_{{\bf m}_{n}}(T_{{\bf
m}_n}(a) \xi)$, $a\in U'_q({\rm so}_{n})$, also belong to these
subspaces ${\cal H}_{{\bf m}^{\pm j}_n}$, respectively.}
 \medskip

{\it Proof.} The assertion follows from the definition of vector
operators and from formula (17).

\bigskip

\centerline{\sc 5. Auxiliary propositions}
\medskip

 \noindent
As stated above, the algebra $U'_q({\rm so}_{n})$ has a
commutative subalgebra ${\cal A}$ generated by the elements
$I_{2s,2s-1}$, $s=1,2,\cdots ,r$, where $r=\lfloor n/2\rfloor$ is
the integral part of $n/2$.
\medskip

{\bf Proposition 2.} (a) {\it If $T$ is a finite dimensional
representation of the algebra $U'_q({\rm so}_n)$, then the
operators
\[
T(I_{21}), T(I_{43}),\cdots ,T(I_{2k,2k-1}),
\]
where $n=2k$ or $n=2k+1$, are simultaneously diagonalizable}.

(b) {\it Possible eigenvalues of any of these operators can be
only as ${\rm i}[m]$, $m\in \frac 12 {\Bbb Z}$, ${\rm
i}=\sqrt{-1}$, or $[m]_+$, $m\in {\Bbb Z}+\frac 12$, where}
\[
[m]\equiv [m]_q=\frac{q^m- q^{-m}}{q-q^{-1}},\ \ \ \
[m]_+=\frac{q^m+ q^{-m}}{q-q^{-1}}.
\]

{\it Proof.} This proposition is true for the algebra $U'_q({\rm
so}_3)$. It follows from complete reducibility of finite
dimensional representations of $U'_q({\rm so}_3)$ (see [12]) and
from the fact that representations of the classical and of the
nonclassical types exhaust all irreducible representations of
$U'_q({\rm so}_3)$ (see [11]). Each of the elements $I_{21},
I_{43},\cdots ,I_{2k,2k-1}$ can be included into some subalgebra
$U'_q({\rm so}_3)$ as one of its generating elements. Therefore,
each of the operators $T(I_{2j,2j-1})$, $j=1,2,\cdots , k$, can be
diagonalized and has eigenvalues indicated in assertion (b) . This
means that these operators are semisimple. Semisimple operators on
a finite dimensional space can be simultaneously diagonalized if
they commute with each other. Proposition is proved.
\medskip

Eigenvalues of the form ${\rm i}[m]$ are called {\it eigenvalues
of the classical type}. Eigenvalues of the form $[m]_+$ are called
{\it eigenvalues of the nonclassical type}.
\medskip

{\it Remark.} In the formulation of Proposition 2 we could take
for the algebra $U'_q({\rm so}_{2k+1})$ the operators $T(I_{32}),
T(I_{54}),\cdots ,T(I_{2k+1,2k})$ instead of $T(I_{21}),
T(I_{43}),\cdots ,T(I_{2k,2k-1})$.
\medskip

In Propositions 3--5 below we suppose that the following
assumption is fulfilled: {\it Each finite dimensional
representation of $U'_q({\rm so}_{n-1})$ is completely reducible
and irreducible finite dimensional representations of $U'_q({\rm
so}_{n-1})$ are exhausted by the irreducible representations of
the classical and nonclassical types described in section 3.} Note
that for $U'_q({\rm so}_{3})$ and $U'_q({\rm so}_{4})$ this
assumption is true (see [10]--[12]).
\medskip

{\bf Proposition 3.} {\it The restriction of any irreducible
finite dimensional representation $T$ of the algebra $U'_q({\rm
so}_{n})$ onto the subalgebra $U'_q({\rm so}_{n-1})$ is completely
reducible representation of $U'_q({\rm so}_{n-1})$ and decomposes
into irreducible representations of this subalgebra which belong
only to the classical type or only to the nonclassical type.}
\medskip

{\it Proof.} The restriction of $T$ to the subalgebra $U'_q({\rm
so}_{n-1})$ is completely reducible due to the assumption. Let
$T{\downarrow}_{U'_q({\rm so}_{n-1})} =\bigoplus _i R_i$, where
$R_i$ are irreducible representations of $U'_q({\rm so}_{n-1})$,
and let ${\cal H}=\bigoplus _i {\cal V}_i$ be the corresponding
decomposition of the space ${\cal H}$ of the representation $T$.
The subspaces ${\cal V}_i$ are invariant with respect to the
operators $T(I_{j,j-1})$, $j=2,3,\cdots , n-1$, corresponding to
the elements of $U'_q({\rm so}_{n-1})$. Only the operator
$T(I_{n,n-1})$ maps vectors of any of the subspaces ${\cal V}_i$
to linear combinations of vectors from other subspaces ${\cal
V}_i$. Since the representation $T$ is irreducible, then acting
repeatedly by $T(I_{n,n-1})$ upon any vector of any subspace
${\cal V}_i$ we obtain linear combinations of vectors from all
other subspaces ${\cal V}_i$. Let some irreducible representation
$R_{i_0}$ of $U'_q({\rm so}_{n-1})$ in the decomposition of $T$
belongs to the classical type. We state that then all other
representations $R_i$ in the decomposition belong to the classical
type. This follows from the following reasoning. We take the
operators $T(I_{n,s})$, $s=1,2,\cdots ,n-1$. It follows from the
commutation relations (4)--(6) for the elements $I_{r,s}$, $r>s$,
given in section 2, that these operators constitute a vector
operator for the subalgebra $U'_q({\rm so}_{n-1})$ (generated by
$I_{21}$, $I_{32},\cdots , I_{n-1,n-2}$) acting on the space
${\cal H}$. Then due to the Wigner--Eckart theorem, the action of
operators $T(I_{n,s})$, $s=1,2,\cdots ,n-1$, on vectors of ${\cal
V}_{i_0}$ gives linear combinations of vectors of subspaces ${\cal
V}_i$ on which only irreducible representations of the classical
type are realized. Repeated application of $T(I_{n,s})$ again
gives representations of the same type. Therefore, in this case,
all representations $R_i$ belong to the classical type. If
$R_{i_0}$ belongs to the nonclassical type, then (by the same
reasoning) all representations $R_{i}$ belong to the nonclassical
type. The proposition is proved.
\medskip

Let us write down the decomposition $T{\downarrow}_{U'_q({\rm
so}_{n-1})} =\bigoplus _i R_i$ from the above proof in the form
$T{\downarrow}_{U'_q({\rm so}_{n-1})} =\bigoplus _{{\bf m}_{n-1}}
d_{{\bf m}_{n-1}} T_{{\bf m}_{n-1}}$ if the decomposition contains
representations of the classical type, where $T_{{\bf m}_{n-1}}$
are irreducible representations of $U'_q({\rm so}_{n-1})$ from
section 3 and $d_{{\bf m}_{n-1}}$ are multiplicities of these
representations. If the decomposition contains irreducible
representations of the nonclassical type, then we write
$T{\downarrow}_{U'_q({\rm so}_{n-1})} =\bigoplus _{\epsilon,{\bf
m}_{n-1}} d_{\epsilon,{\bf m}_{n-1}} T_{\epsilon,{\bf m}_{n-1}}$,
where $T_{\epsilon,{\bf m}_{n-1}}$ are irreducible representations
of the nonclassical type.
\medskip

{\bf Proposition 4.} {\it The action of the operator
$T(I_{n,n-1})$ upon a vector of a subspace, on which the
representation $T_{{\bf m}_{n-1}}$ (the representation
$T_{\epsilon,{\bf m}_{n-1}}$) of $U'_q({\rm so}_{n-1})$ is
realized, gives a linear combination of vectors belonging only to
subspaces of the irreducible representations of $U'_q({\rm
so}_{n-1})$ contained in the decomposition into irreducible
components of the tensor product $T_1\otimes T_{{\bf m}_{n-1}}$
(of the tensor product $T_1\otimes T_{\epsilon,{\bf m}_{n-1}}$),
where $T_1$ is the vector representation of $U'_q({\rm
so}_{n-1})$.}
\medskip

{\it Proof.} The operators $T(I_{n,s})$, $s=1,2,\cdots ,n-1$,
constitute a vector operator for the subalgebra $U'_q({\rm
so}_{n-1})$. Now the proposition follows from the Wigner--Eckart
theorem.
\medskip

{\bf Proposition 5.} {\it Let $T$ be a finite dimensional
irreducible representation of $U'_q({\rm so}_n)$. Then all
operators $T(I_{2i,2i-1})$ from Proposition 2 have eigenvalues
only of the classical type or only of the nonclassical type.}
\medskip

{\it Proof.} The proposition is true for the algebra $U'_q({\rm
so}_4)$. Namely, eigenvalues of $T(I_{21})$ and $T(I_{43})$ of an
irreducible representation $T$ of $U'_q({\rm so}_4)$ are of the
classical type if $T$ is a representation of the classical type
and of the nonclassical type if $T$ is a representation of the
nonclassical type (see [10]). We restrict the representation $T$
of $U'_q({\rm so}_n)$ successively to $U'_q({\rm so}_{n-1})$,
$U'_q({\rm so}_{n-2}), \cdots$, $U'_q({\rm so}_{4})$ and decompose
it into irreducible constituents. (Moreover, the chain of these
subalgebras can be taken in such a way that the last subalgebra
$U'_q({\rm so}_{4})$ contains any two fixed neighbouring operators
from Proposition 2(a).) Applying on the first step Proposition 3
we obtain in the decomposition of $T$ irreducible representations
of $U'_q({\rm so}_{n-1})$ all belonging to the classical type or
all belonging to the nonclassical type. Due to the assumption
before Proposition 3 and Remark at the end of section 3, on each
next step we obtain only irreducible representations of the
classical type or only irreducible representations of the
nonclassical type, described in section 3. Thus, restriction of
$T$ onto any subalgebra $U'_q({\rm so}_{4})$ decomposes into
irreducible representations of $U'_q({\rm so}_{4})$ all belonging
to the classical type or all belonging to the nonclassical type.
Our proposition follows from this assertion. Proposition is
proved.
 \medskip

An irreducible representation $T$ of $U'_q({\rm so}_{n})$ for
which all the operators $T(I_{2i,2i-1})$, $i=1,2,\cdots$, $\lfloor
n/2\rfloor$, have eigenvalues of the classical type (of the
nonclassical type) is called a {\it representation of the
classical type} ({\it of the nonclassical type}). The algebra
$U'_q({\rm so}_{n})$ does not have irreducible finite dimensional
representations of other types. In section 3, irreducible
representations of the classical and of the nonclassical type are given.
But we do not know yet that they exhaust all irreducible
representations of these types. Our aim is to prove that the
irreducible representations of section 3 exhaust all irreducible
finite dimensional representations of $U'_q({\rm so}_{n})$.
\bigskip

\centerline{\sc 6. Reduced matrix elements for the classical type
representations}
\medskip

The theorem on classification of irreducible finite dimensional
representations of the algebra $U'_q({\rm so}_{n})$ will be proved
by means of mathematical induction. Namely, we make an assumption
on irreducible finite dimensional representations of the
subalgebra $U'_q({\rm so}_{n-1})$ (which is true for the
subalgebra $U'_q({\rm so}_{4})$) and then prove that this
assumption is true for the algebra $U'_q({\rm so}_{n})$.
\medskip

{\bf Assumption.} {\it Each finite dimensional representation of
$U'_q({\rm so}_{n-1})$ is completely reducible and irreducible
finite dimensional representations of $U'_q({\rm so}_{n-1})$ are
exhausted by irreducible representations of the classical and
nonclassical types described in section 3.}
\medskip

This assumption is true for the algebras $U'_q({\rm so}_3)$ and
$U'_q({\rm so}_4)$ (see [10] and [11]).

As we know from the previous section, irreducible finite
dimensional representations $T$ of $U'_q({\rm so}_n)$ are divided
into two classes -- irreducible representations of the classical
type and irreducible representations of the nonclassical type. For
deriving the theorem on classification of irreducible
representations belonging to the classical type we need the
results on reduced matrix elements of the tensor operator
$T(I_{n,r})$, $k=1,2,\cdots ,n-1$, for the subalgebra $U'_q({\rm
so}_{n-1})$.

Let $T$ be an irreducible finite dimensional representation of
$U'_q({\rm so}_n)$ belonging to the classical type. According to
our assumption and Proposition 3, this representation decomposes
under the restriction onto the subalgebra $U'_q({\rm so}_{n-1})$
as a direct sum of irreducible representations of the classical
type from section 3. For the space ${\cal H}$ of the
representation $T$ we have
$$
{\cal H}=\bigoplus _{{\bf m}_{n-1},i} {\cal V}_{{\bf m}_{n-1},i},
$$
where ${\cal V}_{{\bf m}_{n-1},i}$ is a linear subspace, on which
the irreducible representation $T_{{\bf m}_{n-1}}$ of $U'_q({\rm
so}_{n-1})$ from section 3 is realized, and $i$ separates multiple
irreducible representations of $U'_q({\rm so}_{n-1})$ in the
decomposition. Let
$$
{\cal V}_{{\bf m}_{n-1}}=\bigoplus _{i} {\cal V}_{{\bf
m}_{n-1},i}.
$$
We take a Gel'fand--Tsetlin basis in each subspace ${\cal V}_{{\bf
m}_{n-1},i}$ and denote these basis vectors by $|{\bf
m}_{n-1},i,\alpha \rangle$, where $\alpha\equiv \alpha_{n-2}$ are
the corresponding Gel'fand--Tsetlin tableaux. Then the subspaces
$$
{\cal V}^\alpha_{{\bf m}_{n-1}}=\bigoplus _i {\Bbb C}|{\bf
m}_{n-1},i,\alpha \rangle
$$
can be defined. We know from Proposition 4 that the operator
$T(I_{n,n-1})$ maps the vector $|{\bf m}_{n-1},i,\alpha \rangle$
into a linear combination of vectors of the subspaces ${\cal
V}_{{\bf m}_{n-1}}$ and ${\cal V}_{{\bf m}^{\pm s}_{n-1}}$,
$s=1,2,\cdots , k$, where $n-1=2k$ or $n-1 =2k+1$. Since the
operator $T(I_{n,n-1})$ commutes with all the operators
$T(I_{s,s-1})$, $s=2,3,\cdots ,n-2$ (that is, with operators
corresponding to elements of the subalgebra $U'_q({\rm
so}_{n-2})$), it maps the subspace ${\cal V}^\alpha_{{\bf
m}_{n-1}}$ into a sum of subspaces ${\cal V}^\alpha_{{\bf
m}'_{n-1}}$ with the same $\alpha$.

Due to Proposition 4 and Wigner--Eckart theorem (see formula
(17)), the action of the operator $T(I_{n,n-1})$ on the subspace
${\cal V}^\alpha_{{\bf m}_{n-1}}$ can be represented in the form
$$
T(I_{2p+2,2p+1})\downarrow _{{\cal V}^\alpha_{{\bf m}_{2p+1}}}
=\sum_{j=1}^p \left( \prod _{r=1}^p
[l_{j,2p+1}+l_{r,2p}][l_{j,2p+1}-l_{r,2p}]\right)^{1/2} \rho_j
({\bf m}_{2p+1}) +
$$  $$
+\sum_{j=1}^p \left( \prod _{r=1}^p
[l_{j,2p+1}+l_{r,2p}-1][l_{j,2p+1}-l_{r,2p}-1]\right) ^{1/2}\tau_j
({\bf m}_{2p+1}) +\left( \prod _{r=1}^p [l_{r,2p}]\right) \sigma
({\bf m}_{2p+1}) \eqno (18)
$$
if $n=2p+2$ and in the form
$$
T(I_{2p+1,2p})\downarrow _{{\cal V}^\alpha_{{\bf m}_{2p}}} =\sum
_{j=1}^p \left( \prod _{r=1}^{p-1}
[l_{j,2p}+l_{r,2p-1}][l_{j,2p}-l_{r,2p-1}+1]\right)^{1/2} \rho'_j
({\bf m}_{2p}) +
$$  $$
+\sum _{j=1}^p \left( \prod _{r=1}^{p-1}
[l_{j,2p}+l_{r,2p-1}-1][l_{j,2p}-l_{r,2p-1}]\right) ^{1/2}
\tau'_j ({\bf m}_{2p})    \eqno (19)
$$
if $n=2p+1$, where $\rho_j ({\bf m}_{2p+1})$,
$\rho'_j ({\bf m}_{2p})$, $\tau_j ({\bf m}_{2p+1})$,
$\tau'_j ({\bf m}_{2p})$ and
$\sigma ({\bf m}_{2p+1})$ are the operators such that
$$
\rho_j ({\bf m}_{2p+1}): {\cal V}^\alpha_{{\bf m}_{2p+1}}\to {\cal
V}^\alpha_{{\bf m}^{+j}_{2p+1}}, \ \ \ \ \rho'_j ({\bf m}_{2p}):
{\cal V}^\alpha_{{\bf m}_{2p}}\to {\cal V}^\alpha_{{\bf
m}^{+j}_{2p}},
$$   $$
\tau_j ({\bf m}_{2p+1}): {\cal V}^\alpha_{{\bf m}_{2p+1}}\to {\cal
V}^\alpha_{{\bf m}^{-j}_{2p+1}},  \ \ \ \ \tau'_j ({\bf m}_{2p}):
{\cal V}^\alpha_{{\bf m}_{2p}}\to {\cal V}^\alpha_{{\bf
m}^{-j}_{2p}},
$$ $$
\sigma ({\bf m}_{2p+1}): {\cal V}^\alpha_{{\bf m}_{2p+1}}\to {\cal
V}^\alpha_{{\bf m}_{2p+1}}
$$
(they are the operators $A_{{\bf m}_{n-1}}^{{\bf m}^{\pm
j}_{n-1}}$ and $A_{{\bf m}_{n-1}}^{{\bf m}_{n-1}}$ from section
4). The last summand in (18) must be omitted if $l_{p,2p+1}=1$ (in
this case the representation $T_{{\bf m}_{2p+1}}$ does not occur
in the tensor product $T_1\otimes T_{{\bf m}_{2p+1}}$). The
coefficients in (18) and (19) are the corresponding
Clebsch--Gordan coefficients of the algebra $U'({\rm so}_{n-1})$
taken from [14]. As we know from the Wigner--Eckart theorem,
$\rho_j ({\bf m}_{2p+1})$, $\rho'_j ({\bf m}_{2p})$, $\tau_j ({\bf
m}_{2p+1})$, $\tau'_j ({\bf m}_{2p})$ and $\sigma ({\bf
m}_{2p+1})$ are independent of $\alpha$. A dependence on $\alpha$
is contained in the Clebsch--Gordan coefficients.

Let us first consider the case of the algebra $U'_q({\rm
so}_{2p+2})$. We act by both parts of the relation
$$
I_{2p+1,2p}I^2_{2p+2,2p+1}-(q+q^{-1})I_{2p+2,2p+1}
I_{2p+1,2p}I_{2p+2,2p+1}
+ I^2_{2p+2,2p+1}I_{2p+1,2p} =-I_{2p+1,2p},
$$
taken for the representation $T$, upon vectors of the subspace
${\cal V}^\alpha_{{\bf m}_{2p+1}}$ with fixed ${\bf m}_{2p+1}$ and
$\alpha$, and take into account formula (18). Comparing terms with
the same resulting subspaces ${\cal V}^\alpha_{{\bf m}'_{2p+1}}$,
we obtain for $\rho_j ({\bf m}_{2p+1})$, $\tau_j ({\bf m}_{2p+1})$
and $\sigma ({\bf m}_{2p+1})$ the relations
$$
[l_{i,2p+1}-l_{j,2p+1}+1]\rho_j ({\bf m}^{+i}_{2p+1})\rho_i ({\bf
m}_{2p+1})- [l_{i,2p+1}-l_{j,2p+1}-1]\rho_i ({\bf
m}^{+j}_{2p+1})\rho_j ({\bf m}_{2p+1})=0, \eqno (20)
$$   $$
[l_{i,2p+1}+l_{j,2p+1}]\tau_i ({\bf m}^{+j}_{2p+1}) \rho_j ({\bf
m}_{2p+1})- [l_{i,2p+1}+l_{j,2p+1}-2] \rho_j ({\bf
m}^{-i}_{2p+1})\tau_i ({\bf m}_{2p+1})=0,  \eqno (21)
$$   $$
[l_{i,2p+1}-l_{j,2p+1}+1]\tau_i ({\bf m}^{-j}_{2p+1}) \tau_j ({\bf
m}_{2p+1})- [l_{i,2p+1}-l_{j,2p+1}-1]\tau_j ({\bf m}^{-i}_{2p+1})
\tau_i ({\bf m}_{2p+1})=0,  \eqno (22)
$$  $$
[l_{j,2p+1}+1]\sigma ({\bf m}^{+j}_{2p+1})\rho_j ({\bf m}_{2p+1})-
[l_{j,2p+1}-1]\rho_j ({\bf m}_{2p+1})  \sigma ({\bf m}_{2p+1})=0,
\eqno (23)
$$  $$
[l_{j,2p+1}] \tau_j ({\bf m}_{2p+1}) \sigma ({\bf m}_{2p+1})-
[l_{j,2p+1}-2] \sigma ({\bf m}^{-j}_{2p+1}) \tau_j ({\bf
m}_{2p+1})=0,  \eqno (24)
$$  $$
\sum_{i=1}^p \biggl( -[2l_{i,2p+1}+1]
 \prod_{r=1\atop r\ne k}^p ([l_{i,2p+1}]^2-[l_{r,2p}]^2)
 \tau_i({\bf m}^{+i}_{2p+1})
 \rho_i ({\bf m}_{2p+1})+
$$  $$
+[2l_{i,2p+1}-3]
 \prod_{r=1\atop r\ne k}^p ([l_{i,2p+1}-1]^2
 -[l_{r,2p}]^2) \rho_i({\bf m}^{-i}_{2p+1})
 \tau_i ({\bf m}_{2p+1}) \biggr) +
\prod_{r=1\atop r\ne k}^p [l_{r,2p}]^2\cdot \sigma^2 ({\bf
m}_{2p+1})=-E, \eqno (25)
$$
where $i\ne j$, $E$ is the unit operator on ${\cal V}^\alpha_{{\bf
m}_{2p+1}}$ and $k$ is a fixed number from the set $\{ 1,2,\cdots
,p\}$. Note that the last term on the left hand side of (25) must
be omitted if $l_{p,2p+1}=1$.

The irreducible representations $T_{{\bf m}_{2p+1}}$ of $U'_q({\rm
so}_{2p+1})$ under restriction to $U'_q({\rm so}_{2p})$ decompose
into irreducible representations $T_{{\bf m}_{2p}}$ of this
subalgebra such that the numbers ${\bf m}_{2p}$ satisfy the
inequalities determined by the Gel'fand--Tsetlin tableaux (see
section 3). Under this, each of the numbers $l_{r,2p}$ runs over a
certain set of values. Assuming that no of $l_{r,2p}$, $r\ne p$,
is a constant for the representation $T_{{\bf m}_{2p+1}}$, we
equate in (25) terms with the same dependence on $[l_{r,2p}]^2$,
$r=1,2,\cdots ,p$, and obtain the relations
$$
\sum _{i=1}^p (-1)^p \Bigl( [2l_{i,2p+1}+1] \tau_i({\bf
m}^{+i}_{2p+1}) \rho_i ({\bf m}_{2p+1})- [2l_{i,2p+1}-3]
\rho_i({\bf m}^{-i}_{2p+1}) \tau_i ({\bf m}_{2p+1})\Bigr)
=-\sigma^2({\bf m}_{2p+1}),
 \eqno (26)
$$  $$
\sum _{i=1}^p \Bigl([2l_{i,2p+1}+1][l_{i,2p+1}]^{2(p-\nu-1)}
\tau_i({\bf m}^{+i}_{2p+1}) \rho_i ({\bf m}_{2p+1})-
$$  $$
-[2l_{i,2p+1}-3][l_{i,2p+1}-1]^{2(p-\nu-1)} \rho_i({\bf
m}^{-i}_{2p+1}) \tau_i ({\bf m}_{2p+1})\Bigr)=0,\ \ \ \nu
=1,2,\cdots ,p-2, \eqno (27)
$$  $$
\sum _{i=1}^p \Bigl([2l_{i,2p+1}+1][l_{i,2p+1}]^{2p-2}
\tau_i({\bf m}^{+i}_{2p+1}) \rho_i ({\bf m}_{2p+1})-
$$  $$
-[2l_{i,2p+1}-3][l_{i,2p+1}-1]^{2p-2} \rho_i({\bf m}^{-i}_{2p+1})
\tau_i ({\bf m}_{2p+1})\Bigr)=E. \eqno (28)
$$
If $s$ parameters $l_{r,2p}$, $r\ne p$, are constant for the
representation $T_{{\bf m}_{2p+1}}$, then the corresponding
$\rho_r ({\bf m}_{2p+1})$ and $\tau_r ({\bf m}_{2p+1})$ vanish and
the number of the relations (27) and (28) is decreased by $s$.

In a similar way it is proved that $\rho'_i ({\bf m}_{2p})$ and
$\tau'_i ({\bf m}_{2p})$ from formula (19) satisfy the relations
$$
[l_{i,2p}-l_{j,2p}+1]\rho'_j ({\bf m}^{+i}_{2p})\rho'_i ({\bf
m}_{2p})- [l_{i,2p}-l_{j,2p}-1]\rho'_i ({\bf m}^{+j}_{2p})\rho'_j
({\bf m}_{2p})=0,\ \ \ i\ne j, \eqno (29)
$$   $$
[l_{i,2p}+l_{j,2p}+1]\tau'_i ({\bf m}^{+j}_{2p}) \rho'_j ({\bf
m}_{2p})- [l_{i,2p}+l_{j,2p}-1] \rho'_j ({\bf m}^{-i}_{2p})\tau'_i
({\bf m}_{2p})=0,\ \ \ i\ne j,  \eqno (30)
$$   $$
[l_{i,2p}-l_{j,2p}+1]\tau'_i ({\bf m}^{-j}_{2p}) \tau'_j ({\bf
m}_{2p})- [l_{i,2p}-l_{j,2p}-1]\tau'_j ({\bf m}^{-i}_{2p}) \tau'_i
({\bf m}_{2p})=0,\ \ \ i\ne j,  \eqno (31)
$$ $$
\sum _i \biggl( -\frac{[2l_{i,2p}+2]}{[l_{i,2p}][l_{i,2p}+1]}
 \prod_{r=1}^{p-1} \Bigl([l_{i,2p}][l_{i,2p}+1]
 -[l_{r,2p-1}][l_{r,2p-1}-1]  \Bigr) \tau'_i({\bf m}^{+i}_{2p})
 \rho'_i ({\bf m}_{2p})+
$$  $$
+\frac{[2l_{i,2p}-2]}{[l_{i,2p}][l_{i,2p}-1]}
 \prod_{r=1}^{p-1} \Bigl([l_{i,2p}][l_{i,2p}-1]-
 [l_{r,2p-1}][l_{r,2p-1}-1] \Bigr) \rho'_i({\bf m}^{-i}_{2p})
 \tau'_i ({\bf m}_{2p}) \biggr) =-E, \eqno (32)
$$
and the last equality leads to the system of equations

$$
\sum _{i=1}^p \Bigl(
[2l_{i,2p}+2]\bigl([l_{i,2p}][l_{i,2p}+1]\bigr)^{p-\nu-2}
\tau'_i({\bf m}^{+i}_{2p}) \rho'_i ({\bf m}_{2p})-
$$  $$
-[2l_{i,2p}-2]\bigl([l_{i,2p}][l_{i,2p}-1]\bigr)^{p-\nu-2}
\rho'_i({\bf m}^{-i}_{2p}) \tau'_i ({\bf m}_{2p})\Bigr) =0, \ \ \
\nu =1,2,\cdots ,p-1, \eqno (33)
$$  $$
\sum _{i=1}^p \Bigl(
[2l_{i,2p}+2]\bigl([l_{i,2p}][l_{i,2p}+1]\bigr)^{p-2} \tau'_i({\bf
m}^{+i}_{2p}) \rho'_i ({\bf m}_{2p})-
$$  $$
-[2l_{i,2p}-2]\bigl([l_{i,2p}][l_{i,2p}-1]\bigr)^{p-2}
\rho'_i({\bf m}^{-i}_{2p}) \tau'_i ({\bf m}_{2p})\Bigr) =E. \eqno
(34)
$$

It follows from the last relations of section 4 that for any $a\in
U'_q({\rm so}_{2p+1})$ the operators $\rho_i ({\bf m}_{2p+1})$,
$\tau_i ({\bf m}_{2p+1})$ and $\sigma({\bf m}_{2p+1})$ satisfy the
relations
$$
T_{{\bf m}_{2p+1}}(a)\sigma ({\bf m}_{2p+1})=\sigma ({\bf
m}_{2p+1}) T_{{\bf m}_{2p+1}}(a),  \eqno (35)
$$  $$
\rho_i ({\bf m}^{-i}_{2p+1})\tau_i ({\bf m}_{2p+1}) T_{{\bf
m}_{2p+1}}(a)= T_{{\bf m}_{2p+1}}(a) \rho_i ({\bf
m}^{-i}_{2p+1})\tau_i ({\bf m}_{2p+1}). \eqno (36)
$$
Similar relations are satisfied by $\rho'_i ({\bf m}_{2p})$ and
$\tau'_i ({\bf m}_{2p})$.

{\it Remark.} Relations (20)--(25) and relations (29)--(32) are
consequences of the relation (2) with $i=n-1$. Other relations
from (1)--(3) containing $I_{n,n-1}$ are satisfied by the
operators (18) and (19). It is a consequence of the fact that
$I_{n,n-1}$ is a component of the vector operator.
\medskip

{\bf Proposition 6.} {\it Let $\xi \in {\cal H}$ belong to a
subspace ${\cal H}_{{\bf m}_{2p+1}}$, on which the irreducible
representation $T_{{\bf m}_{2p+1}}$ of $U'_q({\rm so}_{2p+1})$ is
realized. Then $\rho_j({\bf m}_{2p+1})\xi \in {\cal H}_{{\bf
m}^{+j}_{2p+1}}$ and $\tau_j({\bf m}_{2p+1})\xi \in {\cal H}_{{\bf
m}^{-j}_{2p+1}}$, where ${\cal H}_{{\bf m}^{\pm j}_{2p+1}}$ are
subspaces of ${\cal H}$, on which the irreducible representations
$T_{{\bf m}^{\pm j}_{2p+1}}$ of $U'_q({\rm so}_{2p+1})$ are
realized, respectively. All the vectors $\rho_j({\bf
m}_{2p+1})(T_{{\bf m}_{2p+1}}(a)\xi)$, $a\in U'_q({\rm
so}_{2p+1})$, and all the vectors $\tau_j({\bf m}_{2p+1})(T_{{\bf
m}_{2p+1}}(a)\xi)$, $a\in U'_q({\rm so}_{2p+1})$, belong to these
subspaces ${\cal H}_{{\bf m}^{+j}_{2p+1}}$ and ${\cal H}_{{\bf
m}^{-j}_{2p+1}}$, respectively.}
 \medskip

This proposition is a corollary of Proposition 1.
 \medskip

{\bf Theorem 1.} {\it If the above assumption is true, then the
restriction of an irreducible representation $T$ of $U'_q({\rm
so}_{n})$ to the subalgebra $U'_q({\rm so}_{n-1})$ contains each
irreducible representation of this subalgebra not more than once.}
\medskip

{\it Proof.} We prove the theorem for the algebra $U'_q({\rm
so}_{2p+2})$. For the algebra $U'_q({\rm so}_{2p+1})$ a proof is
the same. Let us consider the decomposition
$$
T{\downarrow}_{U'_q({\rm so}_{2p+1})}=\bigoplus_{{\bf m}_{2p+1}}
d_{{\bf m}_{2p+1}}T_{{\bf m}_{2p+1}}, \eqno (37)
$$
where $d_{{\bf m}_{2p+1}}$ denotes a multiplicity of the
representation $T_{{\bf m}_{2p+1}}$ in the decomposition. The
decomposition ${\cal H}=\bigoplus_{{\bf m}_{2p+1},\alpha} {\cal
V}^\alpha_{{\bf m}_{2p+1}}$ corresponds to the decomposition (37),
where, as in section 4, $\alpha$ numerates elements of the
Gel'fand--Tsetlin basis for the representation $T_{{\bf
m}_{2p+1}}$.  Let $T_{{\bf m}'_{2p+1}}\equiv T_{{\bf m}^{\rm max}
_{2p+1}} $ be a maximal irreducible representation of $U'_q({\rm
so}_{2p+1})$ in the decomposition (37), that is, such that $\rho_j
({\bf m}'_{2p+1})=0$, $j=1,2,\cdots ,p$. Due to the relations
(20)--(22) the operators $\rho_i$ and $\rho_j$, as well as the
operators $\rho_i$ and $\tau_j$, $i\ne j$, and the operators
$\tau_i$ and $\tau_j$, commute (up to a constant) with each other.
For this reason, each of the parameters $l_{i,2p+1}$,
$i=1,2,\cdots ,p$, in the set of the representations $T_{{\bf
m}_{2p+1}}$ from the decomposition (37) runs over some set of
numbers independent of values of other parameters $l_{j,2p+1}$,
$j\ne i$.

We take one of the subspaces ${\cal V}^\alpha_{{\bf m}'_{2p+1}}$,
where ${\bf m}'_{2p+1}\equiv {\bf m}^{\rm max}_{2p+1}$. Its
dimension is equal to the multiplicity $d_{{\bf m}'_{2p+1}}$ of
the representation $T_{{\bf m}'_{2p+1}}$ in the decomposition
(37). Then $\sigma({\bf m}'_{2p+1})$ is an operator on ${\cal
V}^\alpha_{{\bf m}'_{2p+1}}$. Clearly, $\sigma({\bf m}'_{2p+1})$
has at least one eigenvector $\xi_0$ in ${\cal V}^\alpha_{{\bf
m}'_{2p+1}}$. According to (35) all the vectors $T_{{\bf
m}'_{2p+1}}(a)\xi_0$, $a\in U'_q({\rm so}_{2p+1})$, are
eigenvectors of $\sigma({\bf m}'_{2p+1})$. The vectors $T_{{\bf
m}'_{2p+1}}(a)\xi_0$, $a\in U'_q({\rm so}_{2p+1})$, constitute a
subspace ${\cal V}^{\rm ir}_{{\bf m}'_{2p+1}}$, where the
irreducible representation $T_{{\bf m}'_{2p+1}}$ of $U'_q({\rm
so}_{2p+1})$ is realized. Let $\xi_j=\tau_j({\bf
m}'_{2p+1})\xi_0$, $j=1,2,\cdots ,p$. Then $\xi_j\in {\cal
V}^\alpha_{{{\bf m}'}^{-j}_{2p+1}}$ and, due to (21),
$\rho_i({{\bf m}'}^{-j}_{2p+1})=0$ for $i\ne j$. It follows from
(24) that $\xi_j$ is an eigenvector of the operator $\sigma({{\bf
m}'}^{-j}_{2p+1})$. Due to Proposition 6, the vector $T_{{\bf
m}'_{2p+1}}(a)\xi_0$ is mapped by the operator $\tau_j({\bf
m}'_{2p+1})$ into the subspace ${\cal V}^{\rm ir}_{{{\bf
m}'}^{-j}_{2p+1}}$. Hence, the operator $\tau_j({\bf m}'_{2p+1})$
maps ${\cal V}^{\rm ir}_{{\bf m}'_{2p+1}}$ into $\{ 0\}$ or into
the subspace ${\cal V}^{\rm ir}_{{{\bf m}'}^{-j}_{2p+1}}$, on
which the irreducible representation $T_{{{\bf m}'}^{-j}_{2p+1}}$
is realized.

Under a restriction to $U_q'({\rm so}_{2p})$, the representation
$T_{{{\bf m}'}_{2p+1}}$ decomposes into a sum of irreducible
representations $T_{{{\bf m}}_{2p}}$, ${\bf
m}_{2p}=(m_{1,2p},\cdots ,m_{p,2p})$. With the numbers $m_{i,2p}$
we associate numbers $l_{i,2p}$ (see section 3). Suppose that no
of $l_{r,2p}$ is a constant for the representation $T_{{\bf
m}'_{2p+1}}$. We apply both sides of the relations (26)--(28) to
the vector $\xi_0$ and obtain $p$ equations with $p$ unknown
$\rho_i({{\bf m}'}^{-i}_{2p+1}) \tau_i ({\bf m}'_{2p+1})\xi_0$,
$i=1,2,\cdots ,p$. (Note that $\rho_j({\bf m}'_{2p+1})=0$,
$j=1,2,\cdots, p$.) Since $l_{1,2p+1}> l_{2,2p+1}>\cdots
>l_{p,2p+1}$ and $q$ is not a root of unity, the form of
coefficients in (26)--(28) shows that the determinant of this
system is not equal to 0. (In fact, this determinant is
proportional to the Vandermond determinant for $[l_{i,2p+1}]^2$,
$i=1,2,\cdots ,p$.) Solving this system we obtain its (unique)
solution. Since the right hand side of (25) is $-E$, this means
that the vectors $\rho_i({{\bf m}'}^{-i}_{2p+1}) \tau_i ({\bf
m}'_{2p+1})\xi_0$, $i=1,2,\cdots ,p$, are multiple to the vector
$\xi_0$. Since $\tau_i ({\bf m}'_{2p+1})\xi_0=\xi_i$ the vector
$\rho_i({{\bf m}'}^{-i}_{2p+1})\xi_i$ is a multiple to the vector
$\xi_0$. Therefore, due to (36) the operator $\rho_i({{\bf
m}'}^{-i}_{2p+1})$ maps the subspace ${\cal V}^{\rm ir}_{{{\bf
m}'}^{-i}_{2p+1}}$ into $\{ 0\}$ or into ${\cal V}^{\rm ir}_{{{\bf
m}'}_{2p+1}}$. If some of the parameters $l_{r,2p}$ are constant,
then the number of equations (26)--(28) is smaller than $p$. As it
is easy to see, in this case the system of equations also has a
unique solution and the conclusion remains true.

Let $\xi_{j,i}=\tau_j ({{\bf m}'}^{-i}_{2p+1})\xi_i$,
$i=1,2,\cdots ,p$. As above, it is shown that the subspace ${\cal
V}^{\rm ir}_{{{\bf m}'}^{-j,-i}_{2p+1}}$ spanned by the vectors
$T_{{{\bf m}'}^{-j,-i}_{2p+1}}\xi_{j,i}$ is irreducible for
$U'({\rm so}_{2p+1})$ and consists of eigenvectors of the operator
$\sigma({{\bf m}'}^{-j,-i}_{2p+1})$. It is mapped by the operator
$\rho_j({{\bf m}'}^{-j,-i}_{2p+1})$ into $\{ 0\}$ or into ${\cal
V}^{\rm ir}_{{{\bf m}'}^{-i}_{2p+1}}$. Moreover, due to (21), up
to a constant we have
$$
\tau_j ({{\bf m}'}^{-i}_{2p+1})\xi_i=\xi_{j,i}=\tau_i ({{\bf
m}'}^{-j}_{2p+1})\xi_j =\xi_{i,j}. \eqno (38)
$$
Hence, the subspaces constructed by means of the vectors
$\xi_{j,i}$ and $\xi_{i,j}$ coincide. Note that if ${{\bf
m}'}^{-i}_{2p+1}$, ${{\bf m}'}^{-j}_{2p+1}$ and ${{\bf
m}'}^{-j,-i}_{2p+1}$ satisfy the dominance conditions, then the
constant in (38) is not vanishing.

We continue this reasoning further applying successively the
operators $\tau_j$ and $\rho_j$ with appropriate values of the
numbers ${\bf m}_{2p+1}$. Due to the relations (20)--(22) the
operators $\rho_i$ and $\rho_j$, as well as the operators $\rho_i$
and $\tau_j$, $i\ne j$, and the operators $\tau_i$ and $\tau_j$,
commute (up to a constant) with each other. Therefore, as a result
of such continuation, we obtain the set of subspaces ${\cal
V}^{\rm ir}_{{{\bf m}}_{2p+1}}$ of the representation space ${\cal
H}$, on which nonequivalent irreducible representations of the
subalgebra $U'_q({\rm so}_{2p+1})$ are realized and which consist
of eigenvectors of the operators $\sigma({{\bf m}}_{2p+1})$. These
subspaces are mapped by the operators $\rho_i$ and $\tau_i$ into
subspaces of this set. We consider the subspace ${\cal H}'$ of the
space ${\cal H}$ which is a direct sum of these subspaces ${\cal
V}^{\rm ir}_{{{\bf m}}_{2p+1}}$. It follows from the expression
(18) for $T(I_{2p+2,2p+1})$ that this operator leaves ${\cal H}'$
invariant. Due to irreducibility of the representation $T$ we have
${\cal H}'={\cal H}$. This complete a proof for the algebra
$U'_q({\rm so}_{2p+2})$. As is noted above, for $U'_q({\rm
so}_{2p+1})$ a proof is the same. The only difference is that
instead of relations (20)--(28) we have to use relations
(29)--(34). Theorem is proved.
\medskip

The fact that any irreducible representation $T$ of $U'_q({\rm
so}_{n})$ contains each irreducible representation of the
subalgebra $U'_q({\rm so}_{n-1})$ not more than once means that
the operators $\rho_j ({\bf m}_{2p+1})$, $\tau_j ({\bf
m}_{2p+1})$, $\sigma_j ({\bf m}_{2p+1})$, $\rho'_j ({\bf m}_{2p})$
and $\tau'_j ({\bf m}_{2p})$ in (18) and (19) are numerical
functions. Thus, the formula (18) can be represented in the form
$$
T(I_{2p+2,2p+1})| {\bf m}_{2p+1},\alpha\rangle =\sum _j \left(
\prod _{r=1}^p ([l_{j,2p+1}]^2-[l_{r,2p}]^2)\right)^{1/2}\rho_j
({\bf m}_{2p+1}) | {\bf m}^{+j}_{2p+1},\alpha\rangle
$$  $$
+ \sum _j \left( \prod _{r=1}^p
([l_{j,2p+1}-1]^2-[l_{r,2p}]^2\right) ^{1/2}\tau_j ({\bf
m}_{2p+1}) | {\bf m}^{-j}_{2p+1},\alpha\rangle +\left( \prod
_{r=1}^p [l_{r,2p}]\right) \sigma ({\bf m}_{2p+1}) | {\bf
m}^{+j}_{2p+1},\alpha\rangle \eqno (39)
$$
and the formula (19) in the form
$$
T(I_{2p+1,2p}) | {\bf m}_{2p},\alpha\rangle = \sum _j \left( \prod
_{r=1}^{p-1}
([l_{j,2p}+1/2]^2-[l_{r,2p-1}-1/2]^2)\right)^{1/2}\rho'_j ({\bf
m}_{2p})| {\bf m}^{+j}_{2p},\alpha\rangle +
$$  $$
+\sum _j \left( \prod _{r=1}^{p-1}
([l_{j,2p}-1/2]^2-[l_{r,2p-1}-1/2]^2)\right) ^{1/2}\tau'_j ({\bf
m}_{2p})| {\bf m}^{-j}_{2p},\alpha\rangle , \eqno (40)
$$
where $\rho_j ({\bf m}_{2p+1})$, $\tau_j ({\bf m}_{2p+1})$,
$\sigma_j ({\bf m}_{2p+1})$, $\rho'_j ({\bf m}_{2p})$ and $\tau'_j
({\bf m}_{2p+1})$ are appropriate numerical functions.
 \medskip

{\it Remark.} We have seen under proving Theorem 1 that in the set
of the representations $T_{{\bf m}_{2p+1}}$ from the decomposition
(37) each of the parameters $m_{i,2p+1}$, $i=1,2,\cdots,p$, runs
over some set of numbers independent of values of other parameters
$m_{j,2p+1}$, $j\ne i$. It is easy to show by means of formula
(39) that in an irreducible representation $T$ of $U'_q({\rm
so}_{2p+2})$ each $m_{i,2p+1}$, $i=1,2,\cdots,p$, takes all values
from the set $m^{\rm min}_{i,2p+1}, m^{\rm min}_{i,2p+1}+1,\cdots,
m^{\rm max}_{i,2p+1}$ without any omitting. A similar assertion is
true for irreducible finite dimensional representations of
$U'_q({\rm so}_{2p+1})$.
 \medskip

Let us find an explicit form of the functions $\rho_j$, $\tau_j$,
$\sigma$, $\rho'_j$ and $\tau'_j$ from (39) and (40). We first
consider the case of $U'_q({\rm so}_{2p+2})$. From (23) we obtain
the relation $[l_{j,2p+1}+1]\sigma ({\bf
m}^{+j}_{2p+1})=[l_{j,2p+1}-1]\sigma ({\bf m}_{2p+1})$. This means
that $\prod _{j=1}^p [l_{j,2p+1}][l_{j,2p+1}-1]\cdot\sigma ({\bf
m}_{2p+1})$ is independent of $l_{j,2p+1}$, $j=1,2,\cdots ,p$,
that is
$$
\sigma ({\bf m}_{2p+1})= \prod _{j=1}^p
([l_{j,2p+1}][l_{j,2p+1}-1])^{-1}\cdot \sigma, \eqno (41)
$$
where $\sigma$ is a constant. (Note that if $l_{p,2p+1}=1$, then
$\sigma ({\bf m}_{2p+1})\equiv 0$.)

We derive from (20)--(22) the relation
$$
[l_{i,2p+1}-l_{j,2p+1}+1] [l_{i,2p+1}+l_{j,2p+1}+1] \rho_j ({\bf
m}^{+i}_{2p+1})\tau_j ({\bf m}^{+i+j}_{2p+1})=
$$  $$
=[l_{i,2p+1}-l_{j,2p+1}-1] [l_{i,2p+1}+l_{j,2p+1}-1] \rho_j ({\bf
m}_{2p+1})\tau_j ({\bf m}^{+j}_{2p+1}),   \eqno (42)
$$
which shows (after multiplication of both sides by
$[l_{i,2p+1}]^2-[l_{j,2p+1}]^2$) that the expression
$$
([l_{i,2p+1}]^2-[l_{j,2p+1}]^2)([l_{i,2p+1}-1]^2-[l_{j,2p+1}]^2)
\rho_j ({\bf m}_{2p+1})\tau_j ({\bf m}^{+j}_{2p+1})  \eqno (43)
$$
is independent of $l_{i,2p+1}$. Therefore, the expression
$$
\beta_j(l_{j,2p+1})=\rho_j ({\bf m}_{2p+1})\tau_j ({\bf
m}^{+j}_{2p+1})[l_{j,2p+1}]^2[2l_{j,2p+1}-1][2l_{j,2p+1}+1] \times
$$   $$
\times   \prod _{r\ne j}
([l_{r,2p+1}]^2-[l_{j,2p+1}]^2)([l_{r,2p+1}-1]^2-[l_{j,2p+1}]^2)
\eqno (44)
$$
depends only on $l_{j,2p+1}$.

In order to find $\beta_j(l_{j,2p+1})$ we rewrite the relations
(26)--(28) for $\beta_i(l_{i,2p+1})$:
$$
\sum _{i=1}^p \frac{1}{[2l_{i,2p+1}-1]}\left(
\frac{\beta_i(l_{i,2p+1})}{[l_{i,2p+1}]^2 c_i(l_{i,2p+1})} -
\frac{\beta_i(l_{i,2p+1}-1)}{[l_{i,2p+1}-1]^2
c_i(l_{i,2p+1}-1)}\right)=
$$  $$
=(-1)^{p+1}\frac{\sigma^2}{\prod_{r=1}^p
[l_{r,2p+1}]^2[l_{r,2p+1}-1]^2},  \eqno (45)
$$   $$
\sum _{i=1}^p \frac{1}{[2l_{i,2p+1}-1]}\left(
\frac{[l_{i,2p+1}]^{2\nu} \beta_i(l_{i,2p+1})}{ c_i(l_{i,2p+1})} -
\frac{[l_{i,2p+1}-1]^{2\nu}\beta_i(l_{i,2p+1}-1)}{
c_i(l_{i,2p+1}-1)}\right)=0,  \eqno (46)
$$   $$
\nu =0,1,2,\cdots ,p-3,
$$   $$
\sum _{i=1}^p \frac{1}{[2l_{i,2p+1}-1]}\left(
\frac{[l_{i,2p+1}]^{2p-4} \beta_i(l_{i,2p+1})}{ c_i(l_{i,2p+1})} -
\frac{[l_{i,2p+1}-1]^{2p-4}\beta_i(l_{i,2p+1}-1)}{
c_i(l_{i,2p+1}-1)}\right)=1,  \eqno (47)
$$
where
$$
c_i(l_{i,2p+1})=\prod _{r\ne i}
([l_{r,2p+1}]^2-[l_{i,2p+1}]^2)([l_{r,2p+1}-1]^2-[l_{i,2p+1}]^2).
$$

For each fixed $\sigma$, this system of equations has a unique
solution $\beta_i(l_{i,2p+1})$, $i=1,2,\cdots ,p$, since the
determinant of this system is non-vanishing. In order to give this
solution we take into account the constants
$$
l_{r+1,2p+2}=l^{\rm min}_{r,2p+1}-1,\ \ \ \ r=1,2,\cdots ,p,
$$
where $l^{\rm min}_{r,2p+1}$, $r=1,2,\cdots ,p$, are minimal
values of $l_{r,2p+1}$ in the decomposition (37), and represent
$\sigma$ (without loss of a generality) in the form
$$
\sigma ={\rm i}\prod_{r=1}^{p+1} [l_{r,2p+2}] ,  \eqno (48)
$$
where $l_{1,2p+2}$ is a number which is determined by $\sigma$.

From the definition of numbers $l_{r,2p+2}$, $r=2,3,\cdots,p+1$,
and from Remark after Theorem 1 it follows that
$$
l_{2,2p+2}>l_{3,2p+2}>\cdots >l_{p+1,2p+2}.
$$

{\bf Proposition 7.} {\it Solutions of the system (45)--(47) are
given by the expressions
$$
\beta_i(l_{i,2p+1})= \prod_{r=1}^{p+1}
([l_{i,2p+1}]^2-[l_{r,2p+2}]^2)= \sum _{j=0}^{p+1} (-1)^j
e_{p-j+1} ( [l_{1,2p+2}]^2,\cdots ,[l_{p+1,2p+2}]^2)
[l_{i,2p+1}]^{2j}, \eqno (49)
$$
where $e_r(x_1,\cdots x_{p+1})$ are elementary symmetric polynomials in
$x_1,\cdots ,x_{p+1}$.}
\medskip

{\it Proof.} In order to prove this proposition we use the
relations
$$
\sum _{i=1}^s \frac{z_i^m}{\prod_{r=1,r\ne  i}^s (z_i-z_r)}
=\left\{ \matrix{1 & {\rm if} & m=s-1,\cr
        0 & {\rm if} & 0\le m\le s-2,}  \right.
\eqno (50)
$$   $$
\sum _{i=1}^s \frac{1}{z_i\prod_{r=1,r\ne  i}^s (z_i-z_r)}
=\frac{(-1)^{s-1}}{z_1\cdots z_s} \eqno (51)
$$
(see, for example, [25]). We put in these relations $s=2p$ and use
the notations $z_i=x_i, z_{i+p}=y_i$, $i=1,2,\cdots ,p$. Then they
can be written as
$$
\sum _{i=1}^p \frac{1}{x_i{-}y_i}\left( \frac{x_i^m}{\prod_{r\ne
i} (x_r{-}x_i)(y_r{-}x_i)}- \frac{y_i^m}{\prod_{r\ne i}
(x_r{-}y_i)(y_r{-}y_i)} \right)= \left\{ \matrix{1 & {\rm if} &
m=2p-1,\cr
        0 & {\rm if} & 0\le m\le 2p-2,}  \right.
\eqno (52)
$$   $$
\sum _{i=1}^p \frac{1}{x_i-y_i} \left(  \frac{1}{x_i\prod_{r\ne i}
(x_r-x_i)(y_r-x_i)}- \frac{1}{y_i\prod_{r\ne i}
(x_r-y_i)(y_r-y_i)} \right) =\frac{-1}{x_1\cdots x_py_1\cdots
y_p}.  \eqno (53)
$$

We put into the relations (45)--(47) $l_{j,2p+1}=l^{\rm
min}_{j,2p+1}$, $j=1,2,\cdots ,p$, where $l^{\rm min}_{j,2p+1}$ is
a minimal value of $l_{j,2p+1}$ in the decomposition (37). Taking
into account that $\beta_j(l^{\rm min}_{j,2p+1}-1)=0$,
$j=1,2,\cdots ,p$, we see that (45)--(47) turn into a system of
$p$ equations for $\beta_j(l^{\rm min}_{j,2p+1})$, $j=1,2,\cdots
,p$. We substitute into this system the expressions (49) for
$\beta _i(l^{\rm min}_{i,2p+1})$ and then cancel $p-1$ multipliers
from the expression for $\beta _i(l^{\rm min}_{i,2p+1})$ with the
corresponding parts of the expressions for $c_i(l^{\rm
min}_{i,2p+1})$ which are in the denominators. As a result, we
obtain a system of relations which contains only the multiplier
$([l_{1,2p+2}]^2-[l^{\rm min}_{i,2p+1}]^2)$ from $\beta _i(l^{\rm
min}_{i,2p+1})$. Our expressions for $\beta _i(l^{\rm
min}_{i,2p+1})$ are correct if these relations are true. It is
easy to see that they are reduced to the relations (50) and (51)
at $s=p$ if to set $z_i=[l^{\rm min}_{i,2p+1}]^2$, $i=1,2,\cdots
,p$.

Further we prove a correctness of the expressions (49) for $\beta
_i(l_{i,2p+1})$ by induction. Namely, we first put
$l_{j,2p+1}=l^{\rm min}_{j,2p+1}$, $j\ne 1$, and successively
conduct a proof for $\beta _1(l^{\rm min}_{1,2p+1}+1)$, $\beta
_1(l^{\rm min}_{1,2p+1}+2),\cdots$, $\beta _1(l^{\rm
max}_{1,2p+1}-1)$. Then we put $l_{j,2p+1}=l^{\rm min}_{j,2p+1}$,
$j\ne 1,\ 2$, and conduct a proof for $\beta _2(l^{\rm
min}_{2,2p+1}+1)$, $\beta _2(l^{\rm min}_{2,2p+1}+2),\cdots$,
$\beta _2(l^{\rm max}_{2,2p+1}-1)$ under any value of
$l_{1,2p+1}$. We continue this procedure up to $\beta
_p(l_{p,2p+1})$. On each step this proof is conducted by using the
relations (52) and (53). Namely, we put in these relations
$x_i=[l_{i,2p+1}]^2$ and $y_i=[l_{i,2p+1}-1]^2$, then multiply
each of them by the corresponding symmetric polynomial from (49),
and sum up them term-wise in order to obtain the relation (45),
then the relations (46) for $\nu=0,1,2,\cdots , p-3$, and at last
the relation (47). This proves that $\beta_j(l_{j,2p+1})$,
$j=1,2,\cdots ,p$, for given values of $l_{j,2p+1}$ satisfy the
relations (45)--(47). Note that $\beta _i(l^{\rm max}_{i,2p+1})=0$
since in this case $\rho_i({\bf m}^{\rm max}_{2p+1})=0$.
Proposition is proved.
\medskip

Thus, we have found the expressions for $\beta_j(l_{j,2p+1})$,
$j=1,2,\cdots ,p$, depending on $l_{1,2p+2}$, and the
corresponding values of $\sigma$. In order to separate $\rho_j
({\bf m}_{2p+1})$ and $\tau_j ({\bf m}^{+j}_{2p+1})$ in
expression (44) for $\beta_j(l_{j,2p+1})$ we note that these
functions are not determined uniquely by the representation.
Ambiguity in a choice of $\rho_j ({\bf m}_{2p+1})$ and $\tau_j
({\bf m}^{+j}_{2p+1})$ is related to a choice of basis elements.
Namely, in the basis
$$
| {\bf m}_{2p+1},\alpha\rangle '=\prod _{r=1}^p
\omega_r(l_{r,2p+1}) \cdot | {\bf m}_{2p+1},\alpha\rangle ,
$$
where $\omega_r(l_{r,2p+1})$ is a numerical multiplier depending
only on $l_{r,2p+1}$, we obtain somewhat different formulas for
the operator $T(I_{2p+2,2p+1})$. Actually, if to pass to the basis
$\{ | {\bf m}_{2p+1},\alpha\rangle '\}$ in formula (39), then the
coefficient $\sigma ({\bf m}_{2p+1})$ remains without any change,
and $\rho_j ({\bf m}_{2p+1})$ and $\tau_j ({\bf m}_{2p+1})$ are
transformed into
$$
\hat\rho_j ({\bf m}_{2p+1})=\frac{\omega_j(l_{j,2p+1})}
{\omega_j(l_{j,2p+1}+1)} \rho_j ({\bf m}_{2p+1}),\ \ \ \
\hat\tau_j ({\bf m}_{2p+1})=\frac{\omega_j(l_{j,2p+1})}
{\omega_j(l_{j,2p+1}-1)} \tau_j ({\bf m}_{2p+1}).
$$
Moreover, we have
$$
\hat\rho_j ({\bf m}_{2p+1}) \hat\tau_j ({\bf m}^{+j}_{2p+1})
=\rho_j ({\bf m}_{2p+1}) \tau_j ({\bf m}^{+j}_{2p+1}).
$$
It is clear that the multiplier $\omega(l_{j,2p+1})$ can be chosen
in such a way that $\hat\rho_j ({\bf m}_{2p+1})=- \hat\tau_j ({\bf
m}^{+j}_{2p+1})$, that is,
$$
\frac{\omega_j(l_{j,2p+1})} {\omega_j(l_{j,2p+1}+1)} \rho_j ({\bf
m}_{2p+1})=-\frac{\omega_j(l_{j,2p+1}+1)} {\omega_j(l_{j,2p+1})}
\tau_j ({\bf m}^{+j}_{2p+1}).
$$
We obtain from here that
$$
\left( \frac{\omega_j(l_{j,2p+1})} {\omega_j(l_{j,2p+1}+1)}\right)
^2 =-\frac{\tau_j ({\bf m}^{+j}_{2p+1})}{\rho_j ({\bf m}_{2p+1})}.
$$
Taking this relation for $l_{j,2p+1}=l^{\rm min}_{j,2p+1}, l^{\rm
min}_{j,2p+1}+1,l^{\rm min}_{j,2p+1}+2,\cdots$ we find that
$$
\omega_j(l_{j,2p+1})=c\left( \prod_{l=l^{\rm
min}_{j,2p+1}}^{l_{j,2p+1}-1} \frac{\rho_j ({\bf m}_{2p+1})}
{\tau_j ({\bf m}^{+j}_{2p+1})} \right) ^{1/2} ,
$$
where $c$ is a constant. Thus, we may consider that from the very
beginning we have a basis for which
$$
\rho_j ({\bf m}_{2p+1})=-\tau_j ({\bf m}^{+j}_{2p+1}). \eqno (54)
$$
Then it follows from (44), (49) and (54) that
$$
\rho_j ({\bf m}_{2p+1})=\left( \frac{
[l_{j,2p+1}]^{-2}[2l_{j,2p+1}-1]^{-1}
 \prod_{r=1}^{p+1}
([l_{r,2p+2}]^2-[l_{j,2p+1}]^2)}{ [l_{j,2p+1}+1] \prod _{r\ne j}
([l_{r,2p+1}]^2-[l_{j,2p+1}]^2)([l_{r,2p+1}-1]^2-[l_{j,2p+1}]^2)}
\right) ^{1/2}   \eqno (55)
$$
where $l_{r+1,2p+2}=l^{\rm min}_{r,2p+1}-1$, $r=1,2,\cdots p$, and
$l_{1,2p+2}$ is a parameter which together with $l_{r,2p+2}$,
$r=2,3,\cdots, p+1$, must determine irreducible representations.
In the next section we shall find a domain of the parameters
$l_{r,2p+2}$, $r=1,2,\cdots ,p+1$.

Substituting the expressions (54) and (55) for $\rho_j ({\bf
m}_{2p+1})$ and $\tau_j ({\bf m}_{2p+1})$ into (39), we obtain
$$
T(I_{2p+2,2p+1})\vert {\bf m}_{2p+1},\alpha\rangle= \sum^{p}_{j=1}
\frac{B^j_{2p+1}({\bf m}_{2p+1})} {b(l_{j,2p+1})[l_{j,2p+1}]}
\vert {\bf m}^{+j}_{2p+1},\alpha\rangle -
$$
$$
-\sum^{p}_{j=1}\frac {B^j_{2p+1}({\bf m}^{-j}_{2p+1})}
{b(l_{j,2p+1}-1)[l_{j,2p+1}-1]} \vert {\bf m}^{-j}_{2p+1},
\alpha\rangle + {\rm i}\, C_{2p+1}({\bf m}_{2p+1}) \vert {\bf
m}_{2p+1}, \alpha \rangle , \eqno(56)
$$
where $b(l_{j,2p+1})=( [2l_{j,2p+1}+1][2l_{j,2p+1}-1])^{1/2}$ and
$$
B^j_{2p+1}({\bf m}_{2p+1})=\left( \frac{\prod_{i=1}^{p+1}
[l_{i,2p+2}+l_{j,2p+1}] [l_{i,2p+2}-l_{j,2p+1}] \prod_{i=1}^{p}
[l_{i,2p}+l_{j,2p+1}] [l_{i,2p}-l_{j,2p+1}]} {\prod_{i\ne j}^{p}
[l_{i,2p+1}{+}l_{j,2p+1}][l_{i,2p+1}{-}l_{j,2p+1}]
[l_{i,2p+1}{+}l_{j,2p+1}{-}1][l_{i,2p+1}{-}l_{j,2p+1}{-}1]}
\right) ^{1/2} ,
$$
$$
C_{2p+1}({\bf m}_{2p+1}) =\frac{ \prod_{s=1}^{p+1} [ l_{s,2p+2} ]
\prod_{s=1}^{p} [ l_{s,2p} ]} {\prod_{s=1}^{p} [l_{s,2p+1}]
[l_{s,2p+1} - 1] } .
$$
This formula coincides with (9) if to replace $p+1$ by $p$. We
have to determine admissible values of the parameters
$l_{i,2p+2}$, $i=1,2,\cdots, p+1$.

Now we consider the case of $U'_q({\rm so}_{2p+1})$. We have to
find possible expressions for $\rho'_j ({\bf m}_{2p})$ and
$\tau'_j ({\bf m}_{2p})$ in (40).

 We derive from (29)--(31) the relation
$$
[l_{i,2p}+l_{j,2p}] [l_{i,2p}-l_{j,2p}-1] [l_{i,2p}+l_{j,2p}+1]
[l_{i,2p}-l_{j,2p}] \rho'_j ({\bf m}_{2p})\tau'_j ({\bf
m}^{+j}_{2p})=
$$  $$
=[l_{i,2p}+l_{j,2p}] [l_{i,2p}-l_{j,2p}-1] [l_{i,2p}+l_{j,2p}-1]
[l_{i,2p}-l_{j,2p}-2] \rho'_j ({\bf m}^{-i}_{2p})\tau'_j ({\bf
m}^{-i+j}_{2p}),
$$
which shows that the expression
$$
([l_{i,2p}][l_{i,2p}-1]-[l_{j,2p}][l_{j,2p}+1])
([l_{i,2p}+1][l_{i,2p}]-[l_{j,2p}][l_{j,2p}+1]) \rho'_j ({\bf
m}_{2p})\tau'_j ({\bf m}^{+j}_{2p})
$$
is independent of $l_{i,2p}$. Therefore, the expression
$$
\beta'_j(l_{j,2p})=\rho'_j ({\bf m}_{2p})\tau'_j ({\bf
m}^{+j}_{2p})(q^{l_{j,2p}}+q^{-l_{j,2p}})(q^{l_{j,2p}+1}+
q^{-l_{j,2p}-1})
$$   $$
\times \prod_{r\ne j}
([l_{r,2p}][l_{r,2p}-1]-[l_{j,2p}][l_{j,2p}+1])
([l_{r,2p}+1][l_{r,2p}]-[l_{j,2p}][l_{j,2p}+1])
$$
depends only on $l_{j,2p}$. Then we rewrite the relations (33) and
(34) for $\beta'_j(l_{j,2p})$ and in the same way as in
Proposition 7, using the equalities (50) and (52), derive the
following proposition.
\medskip

{\bf Proposition 8.} {\it Solutions of the system of equations for
$\beta'_j(l_{j,2p})$ are given by the expressions
$$
\beta'_j(l_{j,2p})= \prod_{r=1}^{p}
([l_{j,2p}][l_{j,2p}+1]-[l_{r,2p+1}][l_{r,2p+1}-1])=
\prod_{r=1}^{p} [l_{r,2p+1}+l_{j,2p}][l_{r,2p+1}-l_{j,2p}-1]=
$$  $$
=\sum _{j=0}^{p} (-1)^{p-j} e_{p-j} (
[l_{1,2p+1}][l_{1,2p+1}-1],\cdots ,[l_{p,2p+1}][l_{p,2p+1}-1])
([l_{j,2p}][l_{j,2p}+1])^j,
$$
where $l_{i,2p+1}=l^{\rm max}_{i,2p}+1$, $i=1,2,\cdots ,p$, and
$e_r(x_1,\cdots ,x_{p})$ are elementary symmetric polynomials in
$x_1,\cdots ,x_{p}$ .}
\medskip

Separating $\rho'_j ({\bf m}_{2p})$ and $\tau'_j ({\bf
m}^{+j}_{2p})$ from $\beta'_j(l_{j,2p})$ as in the previous case,
for the operator $T(I_{2p+1,2p})$ of an irreducible representation
$T$ of $U'_q({\rm so}_{2p+1})$ we obtain
$$
T(I_{2p+1,2p}) | {\bf m}_{2p},\alpha\rangle = \sum^p_{j=1} \frac{
A^j_{2p}({\bf m}_{2p})} {a(l_{j,2p}) }
            \vert {\bf m}^{+j}_{2p}, \alpha\rangle -
\sum^p_{j=1} \frac{A^j_{2p}({\bf m}^{-j}_{2p})}
 {a(l_{j,2p}-1)} | {\bf m}^{-j}_{2p},\alpha\rangle ,
\eqno (57)
$$
where $a(l_{j,2p})=\{ (q^{l_{j,2p}+1}+q^{-l_{j,2p}-1})
(q^{l_{j,2p}}+q^{-l_{j,2p}})\} ^{1/2}$ and
$$
A^j_{2p}({\bf m}_{2p}) = \left( \frac{\prod_{i=1}^p
[l_{i,2p+1}+l_{j,2p}] [l_{i,2p+1}-l_{j,2p}-1] \prod_{i=1}^{p-1}
[l_{i,2p-1}+l_{j,2p}] [l_{i,2p-1}-l_{j,2p}-1]} {\prod_{i\ne j}^p
[l_{i,2p}+l_{j,2p}][l_{i,2p}-l_{j,2p}]
[l_{i,2p}+l_{j,2p}+1][l_{i,2p}-l_{j,2p}-1]} \right)^{1/2} .
$$

Thus, we derived an explicit form of the operator $T(I_{n,n-1})$
of an irreducible representation of $U'_q({\rm so}_{n})$. In order
to obtain a classification of irreducible representations of the
classical type we have (by using (56) and (57)) to derive a domain
of the parameters $l_{1n},l_{2n},\cdots ,l_{pn}$, $p=\lfloor
n/2\rfloor$.

\bigskip

\centerline{\sc 7. Reduced matrix elements for the nonclassical
type representations}
\medskip

We assume that Assumption of section 6 is acting.
\medskip

{\bf Proposition 9.} {\it Let $T$ be an irreducible finite
dimensional representation of $U'_q({\rm so}_n)$ belonging to the
nonclassical type. Then the decomposition of
$T{\downarrow}_{U'_q({\rm so}_{n-1})}$ into irreducible
constituents contains irreducible representations
$T_{\epsilon,{\bf m}_{n-1}}$ with the same $\epsilon$.}
\medskip

{\it Proof.} The proposition follows from Proposition 4 and from
the fact that the decomposition of the tensor products $T_1\otimes
T _{\epsilon,{\bf m}_{n-1}}$ (where $T_1$ is a vector
representation) into irreducible constituents contains irreducible
representations of the nonclassical type with $\epsilon$
coinciding with $\epsilon$ in $T _{\epsilon,{\bf m}_{n-1}}$.
Proposition is proved.
\medskip

Let $T$ be such as in Proposition 9 and let ${\cal H}$ be a space
on which $T$ acts. Let
$$
{\cal H}=\bigoplus _{{\bf m}_{n-1},i} {\cal V}_{\epsilon,{\bf
m}_{n-1},i}, \eqno (58)
$$
where ${\cal V}_{\epsilon,{\bf m}_{n-1},i}$ is a linear subspace,
on which an irreducible representation $T _{\epsilon,{\bf
m}_{n-1}}$ of $U'_q({\rm so}_{n-1})$ is realized, and $i$
separates multiple irreducible representations in the
decomposition. We also introduce the subspaces
$$
{\cal V}_{\epsilon,{\bf m}_{n-1}}=\bigoplus _{i} {\cal
V}_{\epsilon,{\bf m}_{n-1},i},
$$
We take a Gel'fand--Tsetlin basis in each subspace ${\cal
V}_{\epsilon,{\bf m}_{n-1},i}$ and denote the basis vectors by
$|\epsilon,{\bf m}_{n-1},i,\alpha \rangle$, where $\alpha\equiv
\alpha_{n-2}$ are the corresponding Gel'fand--Tsetlin tableaux.
Let
$$
{\cal V}^\alpha_{\epsilon,{\bf m}_{n-1}}=\bigoplus _i {\Bbb
C}|\epsilon,{\bf m}_{n-1},i,\alpha \rangle .   \eqno (59)
$$
We know from Proposition 4 that the operator $T(I_{n,n-1})$
transforms the vector $|\epsilon,{\bf m}_{n-1},i,\alpha \rangle$
into a linear combination of vectors of the subspaces ${\cal
V}_{\epsilon,{\bf m}_{n-1}}$ and ${\cal V}_{\epsilon,{\bf m}^{\pm
s}_{n-1}}$, $s=1,2,\cdots , k$, where $k=\lfloor \frac 12
(n-1)\rfloor$. Since the operator $T(I_{n,n-1})$ commutes with all
the operators $T(I_{s,s-1})$, $s=2,3,\cdots ,n-2$ (that is, with
operators corresponding to elements of the subalgebra $U'_q({\rm
so}_{n-2})$), it maps subspaces ${\cal V}^\alpha_{\epsilon,{\bf
m}_{n-1}}$ into a sum of subspaces ${\cal V}^\alpha_{\epsilon,{\bf
m}'_{n-1}}$ with the same $\alpha$.

Due to Wigner--Eckart theorem (see formula (17)), the action of
the operator $T(I_{n,n-1})$ on the subspace ${\cal
V}^\alpha_{\epsilon,{\bf m}_{n-1}}$ can be represented in the form
$$
T(I_{2p+2,2p+1}){\downarrow}_{{\cal V}^\alpha_{\epsilon,{\bf
m}_{2p+1}}} =\sum _{j=1}^p \left( \prod _{r=1}^p
[l_{j,2p+1}+l_{r,2p}][l_{j,2p+1}-l_{r,2p}]\right)^{1/2}\rho_j
(\epsilon,{\bf m}_{2p+1}) +
$$  $$
+\sum_{j=1}^p \left( \prod _{r=1}^p
[l_{j,2p+1}+l_{r,2p}-1][l_{j,2p+1}-l_{r,2p}-1]\right) ^{1/2}
\tau_j (\epsilon,{\bf m}_{2p+1}) + \left( \prod _{r=1}^p
[l_{r,2p}]_+\right) \sigma (\epsilon,{\bf m}_{2p+1}), \eqno (60)
$$
if $n=2p+2$ and in the form
$$
T(I_{2p+1,2p}){\downarrow}_{{\cal V}^\alpha_{\epsilon,{\bf
m}_{2p}}} =\sum_{j=1}^p \left( \prod_{r=1}^{p-1}
[l_{j,2p}+l_{r,2p-1}][l_{j,2p}-l_{r,2p-1}+1]\right)^{1/2} \rho'_j
(\epsilon, {\bf m}_{2p}) +
$$  $$
+\sum_{j=1}^p \left( \prod _{r=1}^{p-1}
[l_{j,2p}+l_{r,2p-1}-1][l_{j,2p}-l_{r,2p-1}]\right)^{1/2} \tau'_j
(\epsilon, {\bf m}_{2p})    \eqno (61)
$$
if $n=2p+1$, where $\rho_j (\epsilon,{\bf m}_{2p+1})$, $\rho'_j
(\epsilon,{\bf m}_{2p})$, $\tau_j (\epsilon,{\bf m}_{2p+1})$,
$\tau'_j (\epsilon,{\bf m}_{2p})$ and $\sigma (\epsilon,{\bf
m}_{2p+1})$ are the operators such that
$$
\rho_j (\epsilon,{\bf m}_{2p+1}): {\cal V}^\alpha_{\epsilon,{\bf
m}_{2p+1}}\to {\cal V}^\alpha_{\epsilon,{\bf m}^{+j}_{2p+1}}, \ \
\ \ \rho'_j (\epsilon,{\bf m}_{2p}): {\cal
V}^\alpha_{\epsilon,{\bf m}_{2p}}\to {\cal
V}^\alpha_{\epsilon,{\bf m}^{+j}_{2p}},
$$   $$
\tau_j (\epsilon,{\bf m}_{2p+1}): {\cal V}^\alpha_{\epsilon,{\bf
m}_{2p+1}}\to {\cal V}^\alpha_{\epsilon,{\bf m}^{-j}_{2p+1}},  \ \
\ \
$$$$
\tau'_j (\epsilon,{\bf m}_{2p}): {\cal V}^\alpha_{\epsilon,{\bf
m}_{2p}}\to {\cal V}^\alpha_{\epsilon,{\bf m}^{-j}_{2p}},\quad
\mbox{if $j\ne p$ or $m_{p,2p}\ge \frac32$},
$$$$
\tau'_p (\epsilon,{\bf m}_{2p}): {\cal V}^\alpha_{\epsilon,{\bf
m}_{2p}}\to {\cal V}^\alpha_{\epsilon,{\bf m}_{2p}},\quad \mbox{if
$m_{p,2p}=\frac12$},
$$
$$
\sigma (\epsilon,{\bf m}_{2p+1}): {\cal V}^\alpha_{\epsilon,{\bf
m}_{2p+1}}\to {\cal V}^\alpha_{\epsilon,{\bf m}_{2p+1}}.
$$
The coefficients in (60) and (61) are the corresponding
Clebsch--Gordan coefficients of the algebra $U'({\rm so}_{n-1})$
taken from [14]. As we know from the Wigner--Eckart theorem,
$\rho_j (\epsilon,{\bf m}_{2p+1})$, $\rho'_j (\epsilon,{\bf
m}_{2p})$, $\tau_j (\epsilon,{\bf m}_{2p+1})$, $\tau'_j
(\epsilon,{\bf m}_{2p})$ and $\sigma (\epsilon,{\bf m}_{2p+1})$
are independent of $\alpha$. A dependence on $\alpha$ is contained
in the Clebsch--Gordan coefficients.

We first consider the case of the algebra $U'_q({\rm so}_{2p+2})$.
We act by both parts of the relation
$$
I_{2p+1,2p}I^2_{2p+2,2p+1}-(q+q^{-1})I_{2p+2,2p+1}I_{2p+1,2p}I_{2p+2,2p+1}
+ I^2_{2p+2,2p+1}I_{2p+1,2p} =-I_{2p+1,2p}
$$
upon vectors of the subspace ${\cal V}^\alpha_{\epsilon,{\bf
m}_{2p+1}}$ with fixed $\epsilon$, ${\bf m}_{2p+1}$, $\alpha$ and
take into account formula (60). As a result, we obtain for $\rho_j
(\epsilon,{\bf m}_{2p+1})$, $\tau_j (\epsilon,{\bf m}_{2p+1})$ and
$\sigma (\epsilon,{\bf m}_{2p+1})$ the relations
$$
[l_{i,2p+1}-l_{j,2p+1}+1]\rho_j (\epsilon,{\bf m}^{+i}_{2p+1})
\rho_i (\epsilon,{\bf m}_{2p+1})- [l_{i,2p+1}-l_{j,2p+1}-1] \rho_i
(\epsilon,{\bf m}^{+j}_{2p+1}) \rho_j (\epsilon,{\bf m}_{2p+1})=0,
\eqno (62)
$$   $$
[l_{i,2p+1}+l_{j,2p+1}]\tau_i (\epsilon,{\bf m}^{+j}_{2p+1})
\rho_j (\epsilon,{\bf m}_{2p+1})- [l_{i,2p+1}+l_{j,2p+1}-2] \rho_j
(\epsilon,{\bf m}^{-i}_{2p+1}) \tau_i (\epsilon,{\bf m}_{2p+1})=0,
\eqno (63)
$$   $$
[l_{i,2p+1}-l_{j,2p+1}+1]\tau_i (\epsilon,{\bf m}^{-j}_{2p+1})
\tau_j (\epsilon,{\bf m}_{2p+1})- [l_{i,2p+1}-l_{j,2p+1}-1] \tau_j
(\epsilon,{\bf m}^{-i}_{2p+1}) \tau_i (\epsilon,{\bf m}_{2p+1})=0,
\eqno (64)
$$  $$
[l_{j,2p+1}+1]_+\sigma (\epsilon,{\bf m}^{+j}_{2p+1}) \rho_j
(\epsilon,{\bf m}_{2p+1})- [l_{j,2p+1}-1]_+\rho_j (\epsilon,{\bf
m}_{2p+1}) \sigma (\epsilon,{\bf m}_{2p+1})=0, \eqno (65)
$$  $$
[l_{j,2p+1}]_+ \tau_j (\epsilon,{\bf m}_{2p+1}) \sigma
(\epsilon,{\bf m}_{2p+1})- [l_{j,2p+1}-2]_+ \sigma (\epsilon,{\bf
m}^{-j}_{2p+1}) \tau_j (\epsilon,{\bf m}_{2p+1})=0,  \eqno (66)
$$  $$
\sum _{i=1}^p \biggl( -[2l_{i,2p+1}+1]
 \prod_{r=1\atop r\ne k}^p ([l_{i,2p+1}]_+^2-
 [l_{r,2p}]_+^2)\ \tau_i(\epsilon,{\bf m}^{+i}_{2p+1})
 \rho_i (\epsilon,{\bf m}_{2p+1})+
$$  $$
+[2l_{i,2p+1}-3]
 \prod_{r=1\atop r\ne k}^p ([l_{i,2p+1}-1]_+^2-[l_{r,2p}]_+^2)\
 \rho_i(\epsilon,{\bf m}^{-i}_{2p+1})
 \tau_i (\epsilon,{\bf m}_{2p+1}) \biggr) -
 \prod_{r=1\atop r\ne k}^p [l_{r,2p}]_+^2\cdot \sigma^2
(\epsilon,{\bf m}_{2p+1})=-E, \eqno (67)
$$
where $i\ne j$, $E$ is the unit operator on ${\cal
V}^\alpha_{\epsilon,{\bf m}_{2p+1}}$ and $k$ is a fixed number
from the set $\{ 1,2,\cdots ,p\}$.

The irreducible representations $T_{\epsilon,{\bf m}_{2p+1}}$ of
$U'_q({\rm so}_{2p+1})$ under restriction to $U'_q({\rm so}_{2p})$
decompose into irreducible representations $T_{\epsilon,{\bf
m}_{2p}}$ of this subalgebra such that the numbers ${\bf m}_{2p}$
satisfy the inequalities determined by the Gel'fand--Tsetlin
tableaux. Under this, each of the numbers $l_{r,2p}$ runs over a
certain set of values. Assuming that no of $l_{r,2p}$, $r\ne p$,
is a constant for the representation $T_{\epsilon,{\bf
m}_{2p+1}}$, we equate in (67) terms with the same dependence on
$[l_{r,2p}]^2_+$ and obtain the relations
$$
\sum _{i=1}^p  \Bigl( [2l_{i,2p+1}+1] \tau_i(\epsilon,{\bf
m}^{+i}_{2p+1}) \rho_i (\epsilon,{\bf m}_{2p+1})-
$$$$
-[2l_{i,2p+1}-3] \rho_i(\epsilon,{\bf m}^{-i}_{2p+1}) \tau_i
(\epsilon,{\bf m}_{2p+1})\Bigr)= (-1)^p\sigma^2(\epsilon,{\bf
m}_{2p+1}),
 \eqno (68)
$$  $$
\sum _{i=1}^p \Bigl( [2l_{i,2p+1}+1][l_{i,2p+1}]_+^{2(p-\nu-1)}
\tau_i(\epsilon,{\bf m}^{+i}_{2p+1}) \rho_i (\epsilon,{\bf
m}_{2p+1})-
$$  $$
-[2l_{i,2p+1}-3][l_{i,2p+1}-1]_+^{2(p-\nu-1)} \rho_i(\epsilon,{\bf
m}^{-i}_{2p+1}) \tau_i (\epsilon,{\bf m}_{2p+1})\Bigr) =0,\ \ \
\nu =1,2,\cdots ,p-2, \eqno (69)
$$  $$
\sum _{i=1}^p \Bigl( [2l_{i,2p+1}+1][l_{i,2p+1}]_+^{2p-2}
\tau_i(\epsilon,{\bf m}^{+i}_{2p+1}) \rho_i (\epsilon,{\bf
m}_{2p+1})-
$$  $$
-[2l_{i,2p+1}-3][l_{i,2p+1}-1]_+^{2p-2} \rho_i(\epsilon,{\bf
m}^{-i}_{2p+1}) \tau_i (\epsilon,{\bf m}_{2p+1})\Bigr) =E. \eqno
(70)
$$
If $k$ parameters $l_{r,2p}$, $r\ne p$, are constant for the
representation $T_{\epsilon,{\bf m}_{2p+1}}$, then a number of the
relations (68)--(70) is decreased by $k$.

In a similar way it is proved that $\rho'_i (\epsilon,{\bf
m}_{2p})$ and $\tau'_i (\epsilon,{\bf m}_{2p})$ from formula (61)
satisfy the relations
$$
[l_{i,2p}-l_{j,2p}+1]\rho'_j (\epsilon,{\bf m}^{+i}_{2p}) \rho'_i
(\epsilon,{\bf m}_{2p})- [l_{i,2p}-l_{j,2p}-1] \rho'_i
(\epsilon,{\bf m}^{+j}_{2p})\rho'_j (\epsilon,{\bf m}_{2p})=0,\ \
\ i\ne j, \eqno (71)
$$   $$
[l_{i,2p}+l_{j,2p}+1]\tau'_i (\epsilon,{\bf m}^{+j}_{2p}) \rho'_j
(\epsilon,{\bf m}_{2p})- [l_{i,2p}+l_{j,2p}-1] \rho'_j
(\epsilon,{\bf m}^{-i}_{2p}) \tau'_i (\epsilon,{\bf m}_{2p})=0,\ \
\ i\ne j, \eqno (72)
$$   $$
[l_{i,2p}-l_{j,2p}+1]\tau'_i (\epsilon,{\bf m}^{-j}_{2p}) \tau'_j
(\epsilon,{\bf m}_{2p})- [l_{i,2p}-l_{j,2p}-1] \tau'_j
(\epsilon,{\bf m}^{-i}_{2p}) \tau'_i (\epsilon,{\bf m}_{2p})=0,\ \
\ i\ne j, \eqno (73)
$$ $$
\sum _{i=1}^p \biggl(
-\frac{[2l_{i,2p}+2]}{[l_{i,2p}]_+[l_{i,2p}+1]_+}
 \prod_{r=1}^{p-1} \Bigl([l_{i,2p}]_+[l_{i,2p}+1]_+
 -[l_{r,2p-1}]_+[l_{r,2p-1}-1]_+  \Bigr) \tau'_i(\epsilon,{\bf m}^{+i}_{2p})
 \rho'_i (\epsilon,{\bf m}_{2p})+
$$  $$
+\frac{[2l_{i,2p}-2]}{[l_{i,2p}]_+[l_{i,2p}-1]_+}
 \prod_{r=1}^{p-1} \Bigl([l_{i,2p}]_+[l_{i,2p}-1]_+-
 [l_{r,2p-1}]_+[l_{r,2p-1}-1]_+ \Bigr) \rho'_i(\epsilon,{\bf m}^{-i}_{2p})
 \tau'_i (\epsilon,{\bf m}_{2p}) \biggr) =-E, \eqno (74)
$$
If $l_{p,2p}\equiv m_{p,2p}=1/2$ then $\rho'_p(\epsilon,{\bf
m}^{-p}_{2p})  \tau'_p (\epsilon,{\bf m}_{2p})$ must be replaced
by $(\tau'_p(\epsilon,{\bf m}_{2p}))^2$. Last relation implies the
equalities
$$
\sum _{i=1}^p \Bigl(
[2l_{i,2p}+2]\bigl([l_{i,2p}]_+[l_{i,2p}+1]_+\bigr)^{p-\nu-2}
\tau'_i(\epsilon,{\bf m}^{+i}_{2p}) \rho'_i (\epsilon,{\bf
m}_{2p})+
$$  $$
-[2l_{i,2p}-2]\bigl([l_{i,2p}]_+[l_{i,2p}-1]_+\bigr)^{p-\nu-2}
\rho'_i(\epsilon,{\bf m}^{-i}_{2p}) \tau'_i (\epsilon,{\bf
m}_{2p})\Bigr) =0,\ \ \ \nu =1,2,\cdots ,p-1,
 \eqno (75)
$$  $$
\sum _{i=1}^p \Bigl(
[2l_{i,2p}+2]\bigl([l_{i,2p}]_+[l_{i,2p}+1]_+\bigr)^{p-2}
\tau'_i(\epsilon,{\bf m}^{+i}_{2p}) \rho'_i (\epsilon,{\bf
m}_{2p})-
$$  $$
-[2l_{i,2p}-2]\bigl([l_{i,2p}]_+[l_{i,2p}-1]_+\bigr)^{p-2}
\rho'_i(\epsilon,{\bf m}^{-i}_{2p}) \tau'_i (\epsilon,{\bf
m}_{2p})\Bigr) =E. \eqno (76)
$$

{\bf Theorem 2.} {\it The restriction of a nonclassical type
irreducible representation $T$ of $U'_q({\rm so}_{n})$ to the
subalgebra $U'_q({\rm so}_{n-1})$ contains each irreducible
representation of this subalgebra not more than once.}
\medskip

This theorem is proved (by using relations (62)--(76)) in the same
way as Theorem 1 and we omit this proof.
\medskip

According to this theorem the operators $\rho_j (\epsilon,{\bf
m}_{2p+1})$, $\rho'_j (\epsilon,{\bf m}_{2p})$, $\tau_j
(\epsilon,{\bf m}_{2p+1})$, $\tau'_j (\epsilon,{\bf m}_{2p})$ and
$\sigma (\epsilon,{\bf m}_{2p+1})$ are numerical functions. We
have to find possible expressions for these functions.

First we consider the case of $U'_q({\rm so}_{2p+2})$. We obtain
from (65) that
$$
\sigma ({\bf m}_{2p+1})= \prod _{j=1}^p
([l_{j,2p+1}]_+[l_{j,2p+1}-1]_+)^{-1}\cdot \sigma, \eqno (77)
$$
where $\sigma$ is a constant. As in the case of the
representations of the classical type, from relations (62)--(64)
we derive that the expression
$$
\beta_j(l_{j,2p+1})=\rho_j (\epsilon,{\bf m}_{2p+1})\tau_j
(\epsilon,{\bf
m}^{+j}_{2p+1})[l_{j,2p+1}]_+^2[2l_{j,2p+1}-1][2l_{j,2p+1}+1]
\times
$$   $$
\times   \prod _{r\ne j}
([l_{r,2p+1}]^2-[l_{j,2p+1}]^2)([l_{r,2p+1}-1]^2-[l_{j,2p+1}]^2)
$$
depends only on $l_{j,2p+1}$.

We rewrite the relations (68)--(70) for $\beta_j(l_{j,2p+1})$ and
introduce the notations
$$
l_{r+1,2p+2}=l^{\rm min}_{r,2p+1}-1,\ \ \ \ r=1,2,\cdots ,p.
$$
Then we represent $\sigma$ (without loss of a generality) in the
form
$$
\sigma =\epsilon_{2p+2}\prod_{r=1}^{p+1} [l_{r,2p+2}]_+ ,  \eqno
(78)
$$
where $l_{1,2p+2}$ is a number which is determined by $\sigma$.
\medskip

{\bf Proposition 10.} {\it Solutions of the system of equations
for $\beta_j(l_{j,2p+1})$ are given by the expressions
$$
\beta_i(l_{i,2p+1})= \prod_{r=1}^{p+1}
([l_{i,2p+1}]^2-[l_{r,2p+2}]^2)=\prod_{r=1}^{p+1}
([l_{i,2p+1}]_+^2-[l_{r,2p+2}]_+^2)=
$$   $$
=\sum _{j=0}^{p+1} (-1)^j e_{p-j+1} ( [l_{1,2p+2}]_+^2,\cdots
,[l_{p+1,2p+2}]_+^2) ([l_{j,2p+1}]_+^2)^j,
$$
where $e_r(x_1,\cdots x_{p+1})$ are elementary symmetric
polynomials in $x_1,\cdots ,x_{p+1}$.}
\medskip

This proposition is proved in the same way as Proposition 7 by
using relations (50)--(53).
\medskip

Separation of $\rho_j (\epsilon,{\bf m}_{2p+1})$ and $\tau_j
(\epsilon,{\bf m}^{+j}_{2p+1})$ from $\beta_j(l_{j,2p+1})$ are
fulfilled in the same way as in the case of formula (44) and we
obtain the following formula for $T(I_{2p+2,2p+1})$:
$$
T(I_{2p+2,2p+1})\vert \epsilon,{\bf m}_{2p+1},\alpha\rangle=
\sum^{p}_{j=1} \frac{B^j_{2p+1}({\bf m}_{2p+1})} {b(l_{j,2p+1})
[l_{j,2p+1}]_+} \vert \epsilon,{\bf m}^{+j}_{2p+1},\alpha\rangle -
$$
$$
-\sum^{p}_{j=1}\frac {B^j_{2p+1}(({\bf m}^{-j}_{2p+1})}
{b(l_{j,2p+1}-1)[l_{j,2p+1}-1]_+} \vert \epsilon,{\bf
m}^{-j}_{2p+1}, \alpha\rangle + \epsilon_{2p} {\hat C}_{2p+1}({\bf
m}_{2p+1}) \vert {\bf m}_{2p+1}, \alpha \rangle , \eqno(79)
$$
where $B^j_{2p+1}({\bf m}_{2p+1})$ and $b(l_{j,2p+1})$ are given
by the same expressions as in (56) and
$$
{\hat C}_{2p+1}({\bf m}_{2p+1}) =\frac{ \prod_{s=1}^{p+1} [
l_{s,2p+2} ]_+ \prod_{s=1}^{p} [ l_{s,2p} ]_+} {\prod_{s=1}^{p}
[l_{s,2p+1}]_+ [l_{s,2p+1} - 1]_+ } .
$$
This formula coincides with (15) if to replace $p+1$ by $p$.

Now we consider the case of $U'_q({\rm so}_{2p+1})$. We derive
from the relations (71)--(73) that
$$
\beta'_j(l_{j,2p})=\rho'_j (\epsilon,{\bf m}_{2p})\tau'_j
(\epsilon,{\bf
m}^{+j}_{2p})(q^{l_{j,2p}}-q^{-l_{j,2p}})(q^{l_{j,2p}+1}-
q^{-l_{j,2p}-1})
$$   $$
\times\prod_{r\ne j}
([l_{r,2p}]_+[l_{r,2p}-1]_+-[l_{j,2p}]_+[l_{j,2p}+1]_+)
([l_{r,2p}+1]_+[l_{r,2p}]_+-[l_{j,2p}]_+[l_{j,2p}+1]_+)
$$
depends only on $l_{j,2p}$ (we used here the relation
$[x][x-1]-[y][y-1] =[x]_+[x-1]_+-[y]_+[y-1]_+$). Then we rewrite
the relations (75) and (76) for $\beta'_j(l_{j,2p})$ and, using
the equalities (50) and (52), derive the following proposition.
\medskip

{\bf Proposition 11.} {\it Solutions of the system of equations
for $\beta'_j(l_{j,2p})$ are given by the expressions
$$
\beta'_j(l_{j,2p})= \prod_{r=1}^{p}
([l_{j,2p}]_+[l_{j,2p}+1]_+-[l_{r,2p+1}]_+[l_{r,2p+1}-1]_+),
$$
where $l_{i,2p+1}=l^{\rm max}_{i,2p}+1$, $i=1,2,\cdots ,p$.}
\medskip

We separate $\rho'_j ({\bf m}_{2p})$ and $\tau'_j ({\bf
m}^{+j}_{2p})$ from $\beta'_j(l_{j,2p})$ and obtain for the
operator $T(I_{2p+1,2p})$ of an irreducible representation $T$ of
$U'_q({\rm so}_{2p+1})$ the expression
$$
T(I_{2p+1,2p}) | \epsilon,{\bf m}_{2p},\alpha\rangle =
\delta_{m_{p,2p},1/2}\frac{\epsilon_{2p+1}}{q^{1/2}-q^{-1/2}}
D_{2p}(\alpha_n)|\epsilon,{\bf m}_{2p},\alpha\rangle+
$$  $$
+\sum^p_{j=1} \frac{ A^j_{2p}({\bf m}_{2p})} {a'(l_{j,2p})}
            \vert \epsilon,{\bf m}^{+j}_{2p} \alpha\rangle -
\sum^p_{j=1} \frac{A^j_{2p}({\bf m}^{-j}_{2p})} {a'(l_{j,2p}-1)} |
\epsilon,{\bf m}^{-j}_{2p},\alpha\rangle ,
$$
where $\epsilon_{2p+1}$ takes one of the values $\pm 1$,
$A^j_{2p}({\bf m}_{2p})$ is given by the same expression as in the
case of the formula (57), $a'(l_{j,2p})$ is such as in (14) and
$$
D_{2p}({\bf m}_{2p}) = \frac{\prod_{i=1}^p [l_{i,2p+1}-\frac 12 ]
\prod_{i=1}^{p-1} [l_{i,2p-1}-\frac 12 ] } {\prod_{i=1}^{p-1}
[l_{i,2p}+\frac 12 ] [l_{i,2p}-\frac 12 ] } .
$$

\medskip

\centerline{\sc 8. Complete reducibility}
\medskip

In this section we prove complete reducibility of finite
dimensional representations of $U'_q({\rm so}_{n})$ if Assumption
of section 6 is true. For the algebras $U'_q({\rm so}_{3})$ and
$U'_q({\rm so}_{4})$ this assumption is fulfilled (see [10, 12]).
 \medskip

{\bf Theorem 3.} {\it If Assumption of section 6 is true, then
each finite dimensional representation of $U'_q({\rm so}_{n})$ is
completely reducible.}
\medskip

{\it Proof.} To prove the theorem it is enough to show that every
finite dimensional representation $T$ of $U'_q({\rm so}_{n})$,
containing two irreducible constituents, is completely reducible.
We represent the space ${\cal H}$ of the representation $T$ in the
form ${\cal H}={\cal H} _1\oplus {\cal H}_2$ such that ${\cal
H}_1$ and ${\cal H}_2$ are invariant with respect to $U'_q({\rm
so}_{n-1})$ and on ${\cal H}_1$ and ${\cal H}/{\cal H}_1$
irreducible representations of $U'_q({\rm so}_{n})$ are realized
(we denote them by $T_1$ and $T_2$, respectively). We have to
consider three cases:
 \medskip

{\it Case 1:} One irreducible constituent of $T$ is of the
clasasical type and another of the nonclassical type.

{\it Case 2:} Both irreducible constituents of $T$ are of the
classical type.

{\it Case 3:} Both irreducible constituents of $T$ are of the
nonclassical type.
 \medskip

{\it Proof of case 1.} We restrict the representation $T$ onto
$U'_q({\rm so}_{n-1})$  and decompose it into a direct sum of
irreducible representations of $U'_q({\rm so}_{n-1})$. Then ${\cal
H}$ is the direct sum ${\cal H}={\cal H} _1\oplus {\cal H}_2$,
where ${\cal H} _1$ and ${\cal H}_2$ are sums of the linear
subspaces on which irreducible representations of $U'_q({\rm
so}_{n-1})$ are realized, which belong to the classical type and
to the nonclassical type, respectively. Let $\xi_1\in {\cal H}_1$
transform under an irreducible representation of $U'_q({\rm
so}_{n-1})$. Then due to Proposition 4 and statements of section 4
on decomposition of tensor products of irreducible
representations, $T(I_{n,n-1})\xi_1\in {\cal H}_1$. Similarly, if
$\xi_2\in {\cal H}_2$ transform under an irreducible
representation of $U'_q({\rm so}_{n-1})$, then by the same reason
$T(I_{n,n-1})\xi_2\in {\cal H}_2$. Therefore, ${\cal H}_1$ and
${\cal H} _2$ are invariant (with respect to $U'_q({\rm so}_{n})$)
subspaces of ${\cal H}$. This means that the representation $T$ is
completely irreducible.
 \medskip

{\it Proof of case 2.} Under restriction of the representation $T$
upon $U'_q({\rm so}_{n-1})$, its irreducible constituents $T_1$
and $T_2$ decompose into a direct sum of irreducible
representations of this subalgebra. We denote the corresponding
collections of numbers, characterizing these representations of
$U'_q({\rm so}_{n-1})$, by ${\bf m}_{n-1}$ and $\tilde{\bf
m}_{n-1}$, respectively. The corresponding sets of ${\bf m}_{n-1}$
and of $\tilde{\bf m}_{n-1}$ will be denoted by $\Omega_1$ and
$\Omega_2$, respectively. Since for ${\bf m}_{n-1}\in \Omega_1$
each $m_{i,n-1}$ runs over values independent of values of
$m_{j,n-1}$, $j\ne i$, then in $\Omega_1$ there exists a single
maximal ${\bf m}_{n-1}$ denoted by ${\bf m}^{\rm max}_{n-1}$.
Similarly, in $\Omega_2$ there exists a single $\tilde{\bf m}^{\rm
max}_{n-1}$. We divide case 2 into four subcase:
 \medskip

{\it Subcase 1:} There exists no irreducible representation
$T_{{\bf m}_{n-1}}$ of $U'_q({\rm so}_{n-1})$ with ${\bf
m}_{n-1}\in \Omega_1$ such that $\tilde{\bf m}^{\rm
max}_{n-1}={\bf m}_{n-1}$.

{\it Subcase 2:} The representation $T_{\tilde{\bf m}^{\rm
max}_{n-1}}$ is equivalent to some irreducible representation
$T_{{\bf m}_{n-1}}$, ${\bf m}_{n-1}\in \Omega_1$ and ${\tilde{\bf
m}^{\rm max}_{n-1}}\ne {{\bf m}^{\rm max}_{n-1}}$.

{\it Subcase 3:} ${\tilde{\bf m}^{\rm max}_{n-1}}= {{\bf m}^{\rm
max}_{n-1}}$ and $T_1$ is not equivalent to $T_2$.

{\it Subcase 4:} $T_1$ is equivalent to $T_2$.
 \medskip

We conduct a proof for representations of the algebra $U'_q({\rm
so}_{2p+2})$. For the algebra $U'_q({\rm so}_{2p+1})$ a proof is
similar and we omit it.

Let $\xi$ be a vector of the subspace ${\cal V}^{\rm
irr}_{\tilde{\bf m}^{\rm max}_{2p+1}}$ on which the irreducible
representation $T_{\tilde{\bf m}^{\rm max}_{2p+1}}$ of $U'_q({\rm
so}_{2p+1})$ is realized. A multiplicity of $T_{\tilde{\bf m}^{\rm
max}_{2p+1}}$ in the representation $T{\downarrow}_{ U'_q({\rm
so}_{2p+1})}$ is one. Therefore, $\xi$ is an eigenvector of the
operator $\sigma (\tilde{\bf m}^{\rm max}_{2p+1})$. We the
reasoning of the proof of Theorem 1 acting successively upon $\xi$
by operators $\rho_i$ and $\tau_j$ of section 6 (corresponding to
the appropriate values of $\tilde{\bf m}_{2p+1}$). As a result, we
obtain an invariant (with respect to $U'_q({\rm so}_{2p+2})$)
subspace $\tilde{\cal H}$ of ${\cal H}$ which is a direct sum of
nonequivalent irreducible (with respect to the subalgebra
$U'_q({\rm so}_{2p+1})$) subspaces ${\cal V}^{\rm irr}_{\tilde{\bf
m}_{2p+1}}$. On $\tilde{\cal H}$ the irreducible representation
$T_2$ of $U'_q({\rm so}_{2p+2})$ is realized. Therefore, $T$ is a
direct sum of its subrepresentations $T_1$ and $T_2$.

In subcase 2, ${\tilde{\bf m}^{\rm max}_{2p+1}}$ is not a maximal
set of $(m_{1,2p+1},\cdots,m_{p,2p+1})$ for the representation
$T$. Therefore, there exists $j$, $1\le j\le p$, such that
$\rho_j(\tilde{\bf m}^{\rm max}_{2p+1})\ne 0$. This operator has
one-dimensional kernel ${\cal K}$. We take a vector $\xi\in {\cal
K}$. Thus, $\rho_j({\tilde{\bf m}^{\rm max}_{2p+1}})\xi=0$. Due to
relation (23) $\xi$ is an eigenvector of the operator
$\sigma(\tilde{\bf m}^{\rm max}_{2p+1})$, and due to (20)
$\rho_i(\tilde{\bf m}^{\rm max}_{2p+1})\xi=0$, $1\le i\le p$. Now
a proof is conducted in the same way as in the previous subcase
(by using the reasoning of the proof of Theorem 1).

Since $T_1$ is not equivalent to $T_2$ in subcase 3, we easily
derive from the results of section 6 that for irreducible
representations $T_1$ and $T_2$ the corresponding values
$\sigma({\bf m}^{\rm max}_{2p+1})$ and $\sigma(\tilde{\bf m}^{\rm
max}_{2p+1})$ are different. Therefore, the operator $\sigma({\bf
m}^{\rm max}_{2p+1})$ for the whole representation $T$ is
diagonalizable. We take eigenvectors $\xi_1$ and $\xi_2$ belonging
to different eigenvalues. Then $\rho_j({\bf m}^{\rm
max}_{2p+1})\xi_s=0$, $s=1,2$, for all values of $j$. We act upon
$\xi_1$ and $\xi_2$ by the operators $\rho_i$ and $\tau_j$ and
then, in the same way as in the proof of Theorem 1, obtain two
linear invariant (with respect to $U'_q({\rm so}_{2p+2})$)
subspaces ${\cal H}_1$ and ${\cal H}_2$ of ${\cal H}$ such that
${\cal H}={\cal H}_1\oplus {\cal H}_1$. This proves the theorem
for subcase 3.

For simplicity of notations, in subcase 4 we set
$$
{\bf m}_{2p+1}=(m_{1,2p+1},\cdots ,m_{p,2p+1})\equiv {\bf
m}=(m_{1},\cdots ,m_p),
$$  $$
(l_{1,2p+1},\cdots ,l_{p,2p+1})\equiv (l_{1},\cdots ,l_p).
$$
The operators ${\sigma}({\bf m})$, ${\rho}_j({\bf m})$ and
${\tau}_j({\bf m})$ for the representation $T$ of $U'_q({\rm
so}_{2p+2})$ will be denoted by ${\sigma}^{(T)}({\bf m})$,
${\rho}^{(T)}_j({\bf m})$ and ${\tau}^{(T)}_j({\bf m})$,
respectively. In subcase 4 these operators are of the form
$$
{\sigma}^{(T)}({\bf m}){=} \left(
\begin{array}{ll}
\sigma({\bf m})& \tilde\sigma({\bf m})\\
0 & \sigma({\bf m})
\end{array}
\right),\quad {\rho}^{(T)}_j ({\bf m}){=} \left(
\begin{array}{ll}
\rho_j({\bf m})& \tilde\rho_j({\bf m})\\
0 & \rho_j({\bf m})
\end{array}
\right), \quad {\tau}^{(T)}_j ({\bf m}){=} \left(
\begin{array}{ll}
\tau_j({\bf m})& \tilde\tau_j({\bf m})\\
0 & \tau_j({\bf m})
\end{array}
\right)
\]
where $\sigma({\bf m})$, $\rho_j({\bf m})$, $\tau_j({\bf m})$
$\tilde\sigma({\bf m})$, $\tilde\rho_j({\bf m})$,
$\tilde\tau_j({\bf m})$ are usual functions. Moreover,
$\sigma({\bf m})$, $\rho_j({\bf m})$ and $\tau_j({\bf m})$ are
functions from section 6, corresponding to the irreducible
representation $T_1$. Substituting these expressions for
${\sigma}^{(T)}({\bf m})$ and ${\rho}_j^{(T)}({\bf m})$ into (23),
we obtain identities for elements ${\sigma}({\bf m})$ and
${\rho}_j({\bf m})$, coinciding with (23), and the identities
 $$
[l_{j}+1](\sigma ({\bf m}^{+j})\tilde\rho_j ({\bf m})+
\tilde\sigma ({\bf m}^{+j})\rho_j ({\bf m}))=
[l_{j}-1](\tilde\rho_j ({\bf m}) \sigma ({\bf m})+ \rho_j ({\bf
m})  \tilde\sigma ({\bf m})). \eqno (80)
 $$
The function $\sigma({\bf m})$ corresponds to an irreducible
representation of the algebra $U'_q({\rm so}_{2p+2})$ and is given
by (41) and (48). Using the relation $[l_{j}+1] \sigma ({\bf
m}^{+j})=[l_{j}-1] \sigma ({\bf m})$, following from (23), we
derive from (80) that $[l_{j}+1] \tilde\sigma ({\bf
m}^{+j})=[l_{j}-1] \tilde\sigma ({\bf m})$. Thus, similarly to the
case of $\sigma ({\bf m})$ in section 6 we derive
$$
\tilde\sigma ({\bf m})= \tilde\sigma  \prod _{j=1}^p
([l_{j}][l_{j}-1])^{-1}, \eqno (81)
$$
where $\tilde\sigma$ is a constant. We state that
$\tilde\sigma=0$. In order to show this we remark that if $\tilde
\sigma ({\bf m})=0$ for some ${\bf m}$, then $\tilde\sigma=0$ and
then $\tilde\sigma ({\bf m})=0$ for all ${\bf m}$.

In the case when $l_{p+1,2p+2}=0$, the representation $T_1\sim
T_2$ contains representations of $U'_q({\rm so}_{2p+1})$ with
$l_p=1$. In this case $\sigma=\tilde\sigma=0$.

Let $l_{p+1,2p+2}>0$. It this case $\sigma\ne 0$. From the
relation (25), written for the representation $T$ of $U'_q({\rm
so}_{2p+2})$, we derive that
$$
\sum_{i=1}^p \biggl( -[2l_{i}+1] \prod_{r=1}^{p-1}
([l_{i}]^2-[l_{r,2p}]^2) F_i({\bf m})+
$$  $$
+[2l_{i}-3]
 \prod_{r=1}^{p-1} ([l_{i}-1]^2
 -[l_{r,2p}]^2)  F_i({\bf m}^{-i})   \biggr) +
\prod_{r=1}^{p-1} [l_{r,2p}]^2\cdot 2 \sigma({\bf
m})\tilde\sigma({\bf m})=0,  \eqno (82)
$$
where
$$
F_i({\bf m}):=\tau_i({\bf m}^{+i}) \tilde\rho_i ({\bf
m})+\tilde\tau_i({\bf m}^{+i}) \rho_i ({\bf m}).
$$

Let us consider representations $T_{\bf m}$ of $U'_q({\rm
so}_{2p+1})$ from $T{\downarrow}_{U'_q({\rm so}_{2p+1})}$ with
$m_2,\cdots, m_p$ taking their minimal values. If all $l_{s,2p}$,
$s=1,2,\cdots,p$, are not fixed for these representations, we have
 $$
[2l_1+1] F_1(m_1,m_2^{\rm min},\ldots,m_p^{\rm min}) -[2l_1-3]
F_1(m_1-1,m_2^{\rm min},\ldots,m_p^{\rm min})
$$  $$
+\sum_{i=2}^p [2l_{i}+1] F_i(m_1,m_2^{\rm min},\ldots,m_p^{\rm
min})= (-1)^{p+1} 2 \sigma(m_1,m_2^{\rm min},\ldots,m_p^{\rm
min})\tilde\sigma(m_1,m_2^{\rm min},\ldots,m_p^{\rm min}), \eqno
(83)
$$  $$
[2l_1+1] [l_{1}]^{2\nu} F_1(m_1,m_2^{\rm min},\ldots,m_p^{\rm
min}) -[2l_1-3] [l_{1}-1]^{2\nu} F_1(m_1-1,m_2^{\rm
min},\ldots,m_p^{\rm min})
$$  $$
+\sum_{i=2}^p [2l_{i}+1] [l_{i}]^{2\nu} F_i(m_1,m_2^{\rm
min},\ldots,m_p^{\rm min})=0,\qquad \nu=1,2,\ldots,p-1. \eqno (84)
$$
We sum each equation in (83) and (84)  over $l_1$ from $l_1^{\rm
min}=l_{2,2p+2}+1$  to $l_1^{\rm max}$ with weight coefficients
$[2l_1-1]$ and obtain
$$
\sum_{i=2}^p G_i = 2(-1)^{p+1} \sum_{l_1=l_1^{\rm min}}^{l_1^{\rm
max}} [2l_1-1] \sigma(m_1,m_2^{\rm min},\ldots,m_p^{\rm
min})\tilde\sigma(m_1,m_2^{\rm min},\ldots,m_p^{\rm min}), \eqno
(85)
$$  $$
\sum_{i=2}^p [l_{i}]^{2\nu} G_i=0,\qquad \nu=1,2,\ldots,p-1, \eqno
(86)
$$
where
$$
G_i=\sum_{l_1=l_1^{\rm min}}^{l_1^{\rm max}} [2l_1-1]
 [2l_{i}+1] F_i(m_1,m_2^{\rm min},\ldots,m_p^{\rm min}).
$$
Since the system of homogeneous equations (86) for $G_i$,
$i=2,3,\ldots,p$, has non-vanishing determinant, we get $G_i=0$
and, therefore, (85) gives
$$
\sum_{l_1=l_1^{\rm min}}^{l_1^{\rm max}} [2l_1-1]
\sigma(m_1,m_2^{\rm min},\ldots,m_p^{\rm
min})\tilde\sigma(m_1,m_2^{\rm min},\ldots,m_p^{\rm min})=0.
 $$
Taking into account (41) and (81) we get
\[
0=\sigma\tilde\sigma\sum_{l_1=l_1^{\rm min}}^{l_1^{\rm max}}
\frac{[2l_1-1]}{[l_1]^2[l_1-1]^2}= \sigma
\tilde\sigma\sum_{l_1=l_1^{\rm min}}^{l_1^{\rm max}}
\left(\frac{1}{[l_1-1]^2}-\frac{1}{[l_1]^2}\right)= \sigma
\tilde\sigma \left(\frac{1}{[l_{2,2p+2}]^2}-\frac{1}{[l_1^{\rm
max}]^2}\right).
\]
Since $[l_1^{\rm max}]^2\ne [l_{2,2p+2}]^2$ and $\sigma\ne 0$, we
obtain $\bar\sigma=0$.

If values of $l_{s,2p}$ are fixed in the considered
representations of $U'_q({\rm so}_{2p+1})$, then the number of
relations which follow from (82) and the number of $G_i$ are
decreased by the number of fixed $l_{s,2p}$. Thus, as before, we
get $G_i=0$, $i=2,3,\ldots,p$ and, therefore, $\tilde\sigma=0$.

We have proved that $\tilde\sigma ({\bf m})=0$ for all irreducible
representations $T_{\bf m}$ of $U'_q({\rm so}_{2p+1})$, contained
in the representation $T{\downarrow}_{U'_q({\rm so}_{n})}$. This
means that all operators $\sigma^{(T)} ({\bf m})$ are diagonal and
the further proof of complete reducibility are conducted in the
same way as in the previous subcase.

The case 3 is proved in the same way as the case 2 and we omit
this proof. The theorem is proved.
\medskip

{\bf Corollary.} {\it If irreducible finite dimensional
representations of $U'_q({\rm so}_{n-1})$ are exhausted by
irreducible representations of section 3, then each finite
dimensional representation of $U'_q({\rm so}_{n})$ is completely
reducible.}

\bigskip

\centerline{\sc 9. Classification theorems}
 \medskip

Suppose that Assumption of section 6 is acting.
 \medskip

{\bf Proposition 12.} {\it If Assumption of section 6 is true,
then irreducible finite dimensional representations $T$ of
$U'_q({\rm so}_{n})$ such that the restriction
$T{\downarrow}_{U'_q({\rm so}_{n-1})}$ contains in the
decomposition into irreducible components only representations of
the classical type of $U'_q({\rm so}_{n-1})$ are exhausted by the
representations of the classical type from section 3.}
\medskip

{\it Proof.} We prove the proposition when $n=2p+2$. For $n=2p+1$
a proof is similar.

Let $T$ be a representation of $U'_q({\rm so}_{2p+2})$ from the
formulation of the proposition. Then the functions
$\beta_j(l_{i,2p+1})$, defined by the formula (44), are given by
(49). It was shown above that $T{\downarrow}_{U'_q({\rm
so}_{2p+1})}=\bigoplus _{{\bf m}_{2p+1}} T_{{\bf m}_{2p+1}}$ and
in this decomposition each $m_{r,2p+1}$ runs over the values
$m^{\rm min}_{r,2p+1}, m^{\rm min}_{r,2p+1}+1,\cdots ,m^{\rm
max}_{r,2p+1}$, where $l^{\rm min}_{r,2p+1}=l_{r+1,2p+2}+1$. Due
to properties of the functions $\rho_j$, $\beta_r(l^{\rm
min}_{r,2p+1}+s)\ne 0$ for $s=0,1,\cdots ,l^{\rm
max}_{r,2p+1}-l^{\rm min}_{r,2p+1}-1$ and $\beta_r(l^{\rm
max}_{r,2p+1})=0$. Then it follows from (49) that $l^{\rm
max}_{r,2p+1}=l_{r,2p+2}$, $r\ne 1$. Since $\beta_r(l^{\rm
max}_{1,2p+1})=0$, we find from (49) that $l^{\rm max}_{1,2p+1}$
coincides with $l_{1,2p+2}$ or with $-l_{1,2p+2}$. Therefore,
$l_{1,2p+2}$ is an integer (a half-integer) if $l_{i,2p+2}$,
$i=2,3,\cdots ,p+1$, are integers (half-integers). Moreover,
$l_{1,2p+2}$ may be positive or negative. We see that the formula
for the operator $T(I_{2p+2,2p+1})$ does not change if we replace
$l_{1,2p+2}$ and $l_{p+1,2p+2}$ by $-l_{1,2p+2}$ and
$-l_{p+1,2p+2}$, respectively. Therefore, we may consider that
$l_{1,2p+2}$ is positive and $l_{p+1,2p+2}$ takes positive and
negative values. Now taking into account admissible values for
$l_{i,2p+2}$, $i=1,2,\cdots ,p+1$, and formula (56) for
$T(I_{2p+2,2p+1})$ we see that the representation $T$ coincides
with one of the irreducible representations of the classical type
from section 3.

In order to prove the proposition for representations of the
algebra $U'_q({\rm so}_{2p+1})$ we use the formula of Proposition
8 and formula (57) instead of formulas (49) and (56). Proposition
is proved.
 \medskip

{\bf Proposition 13.} {\it If Assumption of section 6 is true,
then irreducible finite dimensional representations $T$ of
$U'_q({\rm so}_{n})$ such that the restriction
$T{\downarrow}_{U'_q({\rm so}_{n-1})}$ contains in the
decomposition into irreducible components only representations of
the nonclassical type of $U'_q({\rm so}_{n-1})$ are exhausted by
the representations of the nonclassical type of section 3.}
\medskip

Proof of this proposition is the same as that of Proposition 12.
 \medskip

{\bf Theorem 4.} {\it Irreducible finite dimensional
representations of the algebra $U'_q({\rm so}_{n})$ are exhausted
by representations of the classical type and of the nonclassical
type from section 3.}
 \medskip

{\it Proof.} For the algebra $U'_q({\rm so}_{n-1})\equiv U'_q({\rm
so}_{4})$, Assumption of section 6 is true (see [10]). Now the
theorem is easily proved by induction taking into account Theorem
3 and Propositions 12 and 13. Theorem is proved.
 \medskip

{\bf Corollary.} {\it Each finite dimensional representation of
$U'_q({\rm so}_{n})$ is completely reducible.}
 \medskip

{\it Proof.} This assertion follows from Corollary of section 8
and from Theorem 4.

 \bigskip

\noindent {\bf Acknowledgement}

\medskip

The research of N.Z.I. was partially supported by the INTAS grant
No.~03-51-3350.

\end{document}